\documentclass[letterpaper, 10 pt, journal, twoside]{IEEEtran}

\IEEEoverridecommandlockouts                            
\usepackage{cite}
\usepackage{amsmath,amssymb,amsfonts}
\usepackage{graphicx}
\usepackage{textcomp}
\usepackage{dsfont}
\usepackage{color}
\usepackage{epstopdf}        
\usepackage{lscape}
\usepackage{multicol}
\usepackage{multirow}
\usepackage{calligra}
\usepackage{booktabs}
\usepackage{mathtools}
\usepackage{framed} 
\usepackage{empheq}
\usepackage{graphics} % for pdf, bitmapped graphics files
\usepackage{epsfig} % for postscript graphics files
\usepackage{times} % assumes new font selection scheme installed
\usepackage{subfig}
\usepackage{epstopdf}
\usepackage{enumerate}
\usepackage{graphicx} % grÂficos
\usepackage{xcolor}
\usepackage{comment} 
\usepackage{cite}
\usepackage{amsmath,amssymb,amsfonts}
\usepackage{algorithmic}
\usepackage{textcomp}
\usepackage{dsfont}
\usepackage{color}
\usepackage{epstopdf}        
\usepackage{lscape}
\usepackage{multicol}
\usepackage{multirow}
\usepackage{calligra}
\usepackage{booktabs}
\usepackage{mathtools}
\usepackage{framed} 
\usepackage{empheq}
\usepackage{graphics} % for pdf, bitmapped graphics files
\usepackage{epsfig} % for postscript graphics files
\usepackage{times} % assumes new font selection scheme installed
\usepackage{amsmath} % assumes amsmath package installed
\usepackage{amssymb}  % assumes amsmath package installed
\usepackage{amsfonts}  % assumes amsmath package installed
\usepackage{subfig}
\usepackage{epstopdf}
\usepackage{graphicx} % grÂficos
\usepackage{xcolor}
\usepackage{comment} 
\usepackage{algorithm}

\usepackage{bm}

\newcommand{\tcb}{\textcolor{black}}

\newtheorem{theorem}{Theorem}

\newtheorem{lemma}{Lemma}
\newtheorem{corollary}{Corollary}
\newtheorem{definition}{Definition}
\newtheorem{assumption}{Assumption}
\newtheorem{remark}{Remark}
\newtheorem{example}{Example}

\newcommand{\bx}{\mathbf{x}}

\newcommand{\by}{\mathbf{y}}
\newcommand{\bz}{\mathbf{z}}
\newcommand{\bq}{\mathbf{q}}
\newcommand{\bp}{\mathbf{p}}
\newcommand{\be}{\mathbf{e}}
\newcommand{\bbf}{\mathbf{f}}
\newcommand{\bs}{\mathbf{s}}
\newcommand{\bg}{\mathbf{g}}

\newcommand{\bv}{\mathbf{v}}

\newcommand{\bd}{\mathbf{d}}

\newcommand{\bh}{\mathbf{h}}
\newcommand{\bla}{\bm{\lambda}}

\newcommand{\R}{\mathbb{R}}

\newcommand{\sX}{\mathcal{X}}

\newcommand{\sY}{\mathcal{Y}}
\newcommand{\sP}{\mathcal{P}}

\newcommand{\sB}{\mathcal{B}}
\newcommand{\sO}{\mathcal{O}}

\newcommand{\sF}{\mathcal{F}}

\newcommand{\B}{\mathbb{B}}

\newcommand{\bxi}{\bm{\xi}}
\newcommand{\bmu}{\bm{\mu}}
\newcommand{\bth}{\bm{\theta}}

\newcommand*{\QEDB}{\hfill\ensuremath{\square}}%

\begin{document}

\hyphenation{op-tical net-works semi-conduc-tor}

\allowdisplaybreaks

\title{Continuous-Time Zeroth-Order Dynamics with Projection Maps: Model-Free Feedback Optimization with Safety Guarantees}

\author{Xin Chen, Jorge I. Poveda,  Na Li 
\thanks{X. Chen and J. I. Poveda contributed equally. X. Chen is with the Department of Electrical and Computer Engineering, Texas A\&M University, TX, USA. E-mail: {\tt xin\_chen@tamu.edu}. J. I. Poveda is with the Electrical and Computer Engineering Department, University of California, San Diego, CA, USA. E-mail: {\tt poveda@ucsd.edu}. N. Li is with the School of Engineering and Applied Sciences at Harvard University, MA, USA. E-mail: {\tt nali@seas.harvard.edu}. (Corresponding author: Xin Chen, {\tt xin\_chen@tamu.edu}).}
 \thanks{This work was supported in part by NSF CAREER: ECCS 2305756, NSF CMII-2228791, NSF AI Institute: 2112085,
NSF ASCENT: 2328241, NSF CNS: 2003111, and AFOSR grant FA9550-22-1-0211.}
%\thanks{Earlier, partial results of this paper appeared in the Proceedings of the IEEE CDC 2021 \cite{chen2021safe}.
}

\maketitle

\begin{abstract}
This paper introduces a class of model-free feedback methods for solving generic \emph{constrained} optimization problems where the mathematical forms of the cost and constraint functions are not available. The proposed methods, termed Projected Zeroth-Order (P-ZO) dynamics, incorporate \emph{projection maps} into a class of continuous-time zeroth-order dynamics that use direct measurements of the cost function and periodic dithering for the purpose of gradient learning. In particular, the proposed P-ZO algorithms can be interpreted as new extremum-seeking algorithms that autonomously drive an unknown system toward a neighborhood of the set of solutions of an optimization problem using only output feedback, while simultaneously guaranteeing that the input trajectories remain in a feasible set for all times. In this way, the P-ZO algorithms can properly handle hard and asymptotic constraints in model-free optimization problems without using penalty terms or barrier functions. Moreover, the proposed dynamics have suitable robustness properties with respect to small bounded additive disturbances on the states and dynamics, a property that is fundamental for practical real-world implementations. % Since some of the proposed algorithms incorporate discontinuous projection maps, our stability results rely on multi-time scale techniques for well-posed autonomous differential inclusions with fast and oscillating states. 
Additional tracking results for time-varying and switching cost functions are also derived under stronger convexity and smoothness assumptions and using tools from hybrid dynamical systems. %For the smooth P-ZDs, we further extend our results to decentralized settings suitable for cooperative multi-agent optimization problems where each agent only needs to perform local computations and actions. 
Numerical examples are presented throughout the paper to illustrate the above results. 
\end{abstract}

%%%%%%%%%%%%%%%%%%
%%%%%%%%%%%%%%%%%%%
%

% This paper introduces a class of model-free feedback methods for solving generic constrained optimization problems where the mathematical forms of the cost and constraint functions are not available. The proposed methods, termed Projected Zeroth-Order (P-ZO) dynamics, incorporate projection maps into a class of continuous-time zeroth-order dynamics that use direct measurements of the cost function and periodic dithering for the purpose of gradient learning. In particular, the proposed P-ZO algorithms can be interpreted as new extremum-seeking algorithms that autonomously drive an unknown system toward a neighborhood of the set of solutions of an optimization problem using only output feedback, while systematically guaranteeing that the input trajectories remain in a feasible set for all times. In this way, the P-ZO algorithms can properly handle hard and asymptotical constraints in model-free optimization problems without using penalty terms or barrier functions. Moreover, the proposed dynamics have suitable robustness properties with respect to small bounded additive disturbances on the states and dynamics, a property that is fundamental for practical real-world implementations. Additional tracking results for time-varying and switching cost functions are also derived under stronger convexity and smoothness assumptions and using tools from hybrid dynamical systems. Numerical examples are presented throughout the paper to illustrate the above results. 

\section{INTRODUCTION}
\label{sec:introduction}
%%%%%%%%%%%%%%%%
%
\IEEEPARstart{T}{his} paper studies the design of model-free feedback control algorithms for autonomously steering a plant toward the set of solutions of an optimization problem using high-frequency dither signals. %It relates to the paradigm called feedback optimization \cite{he2022model,hauswirth2021optimization} that interconnects optimization iterations in closed-loop with the optimal control of physical plants. 
This type of feedback control design % dates back to the early 1920s for electromechanical systems \cite{Leblanc} and 
has recently attracted considerable attention, due to successful applications in power grids,
 %\cite{hauswirth2021projected}, 
communication networks, %semi-conductor manufacturing systems \cite{ren2012laser}, and
and mobile robots; see \cite{hauswirth2021projected} and references therein. The design of these controllers for practical applications is particularly challenging because of two major obstacles: one is the \emph{lack of accurate models of the system}, as many real-world systems are too complex to derive tractable mathematical equations that accurately describe their behavior in unknown or dynamic environments; the other obstacle is  to meet \emph{safety requirements} by  properly {handling various constraints}, including physical laws, control saturation, capacity and budget limits, etc. This paper introduces a class of algorithms that can overcome both of these obstacles and are suitable for the solution of model-free \emph{constrained} optimization problems describing safety-critical applications.

\vspace{-0.2cm}\noindent 
\subsection{Literature Review} 
To address the problem of unknown system models, real-time \emph{model-free} control and optimization schemes have been extensively studied. In these approaches,
instead of pre-establishing a complex and often \textit{static/stationary} system model from first principles and historical data, adaptive algorithms are used to probe the unknown plant and learn its optimal operation points using real-time output feedback. Such techniques, called \emph{extremum seeking (ES) controllers}, leverage multi-time scale principles to steer dynamical systems to optimal steady state operating points, while preserving closed-loop stability guarantees. ES techniques date back to the early 1920s \cite{Leblanc}. However, the first general stability analysis for nonlinear systems was presented in the 2000s in \cite{Krstic2000} using averaging-based methods, and in \cite{TeelPopovic2001} using sampled-data approaches based on finite-differences approximations. Since these methods rely solely on measurements of the objective function, ES is closely related to discrete-time zeroth-order optimization dynamics \cite{nesterov2017random, chen2022improve}. In the continuous-time domain, ES algorithms have been further advanced during the last two decades using more general analytical and design techniques for ordinary differential equations (ODEs), see \cite{tan2006non,Nesic_ChineseCDC,DurrLieBracket,scheinker2017model,PovedaKrsticFXT}.
 
%In the context of averaging-based ES, the basic idea behind the approach is to multiply the output of the plant with a small periodic dither signal that is also added to the input. When the frequency of the dither is high, and there is enough time scale separation between the plant and the controller, the resulting signal of the product can be seen -on average- as an estimate of the gradient of the objective function. This signal can then be used to close the loop and asymptotically steer the input using an optimization algorithm with integral (I) (or proportional-integral (PI) \cite{guay2016perturbation}) action. %In this way, ES feedback exhibits a close connection with  zeroth-order optimization (ZO) approaches \cite{nesterov2017random,chen2022improve}, traditionally used in the optimization literature for problems where only function evaluations or measurements are available for the user. %Thanks to its model-free feedback nature with closed-loop stability guarantees, ES has been applied to control various black-box or hard-to-model systems, such as voltage control in power systems \cite{chen2021safe}, combustion control of thermal engines \cite{killingsworth2009hcci},  photovoltaic  maximum power point tracking \cite{6362193}, source seeking in mobile robots \cite{Poveda:20TAC,Khong:14}, etc.  
   
However, despite the theoretical advances and practical applications in ES, one of the major challenges of existing schemes is how to guarantee the systematic satisfaction of \emph{hard} and \emph{asymptotic} constraints simultaneously. Hard constraints refer to physical or safety-critical constraints that need to be satisfied by the actions of the controller \emph{at all times}, e.g. saturation or actuator capacity limits, the generation capacity of a power plant, etc. On the other hand, asymptotic constraints refer to soft physical limits or performance requirements %\lina{the word ``artificial'' sounds ``negative'' and miss-leading. I just deleted it} 
  that can be violated temporarily during transient processes but should be met
in the long-term steady state, e.g., the thermal limits of power lines and voltage limits imposed by industrial standards,  the comfortable temperature ranges required in building climate control, etc. 
Properly handling these two types of constraints is essential to ensure stability and optimality in real-time optimization algorithms.

In the context of ES, most of the approaches and stability results have been developed for unconstrained optimization problems. For optimization problems with hard constraints, the majority of the results and schemes have been limited to methods that integrate barrier or penalty functions in the cost \cite{dehaan2005extremum,guay2015constrained,guay2018distributed,tan2013extremum,hazeleger2022sampled,Poveda:20TAC}, which can limit the type and number of constraints that can be handled by the algorithms. In \cite{durr2014extremum,sync_ES_geodesic,Poveda:15,Ochoa2022,PovedaTAC17B}, ES algorithms were introduced to solve optimization problems with constraints defined by certain Euclidean smooth manifolds and Lie Groups. These schemes, however, do not incorporate soft constraints in the optimization problem, and can only handle boundaryless manifolds. Anti-windup techniques in ES for problems that involve saturation were studied in \cite{ESAntiwindup2}, and ES with output constraints were studied in \cite{liao} using boundary tracing techniques. Switching ES algorithms that emulate sliding-mode techniques were also presented in \cite{galarza2022sliding} to handle hard constraints in time-varying problems. More recently, an innovative approach that combines safety filters and ES was introduced in \cite{safeES} using control barrier functions and quadratic programming. To handle soft constraints, ES approaches based on saddle flows have also been studied in \cite{ye2016distributed,durr2013saddle,wang2019distributed,PovedaNaLi2019}. Finally, more closely related to our setting are the works \cite{mills2014constrained,Frihauf12a}, which considered ES algorithms with certain projection maps for scalar problems \cite{mills2014constrained}, and numerical studies of Nash-seeking problems with box constraints \cite[Section V-B]{Frihauf12a}.

\vspace{-0.1cm}
\subsection{Contributions and Organization} 
This paper introduces a class of continuous-time projected zeroth-order (P-ZO) algorithms for solving generic constrained optimization problems with both hard and asymptotic constraints. Based on ES and two different types of \emph{projection maps},
the proposed P-ZO methods can be interpreted as model-free feedback controllers that steer a plant towards the set of solutions of an optimization problem with hard and soft constraints, using only measurements or evaluations of the objective and constraint functions. We explain the main advantages and innovations of the proposed algorithms below:

\vspace{0.1cm}
\noindent 
(a) \emph{Model-Free Methods:} We study a class of optimization algorithms that use only measurements or evaluations of the objective function and the constraints, i.e., zeroth-order (ZO) information. In this way, the algorithms do not require knowledge of the mathematical forms of the expressions that define the optimization problem, or their gradients. We show that, under suitable tuning of the control parameters, the trajectories of the proposed model-free ZO algorithms can approximate the behavior of smooth and non-smooth first-order continuous-time model-based dynamics \cite{nagurney2012projected,gao2003exponential,hauswirth2021projected}. By using real-time output feedback, the proposed algorithms are inherently robust to unknown disturbances. They are also effective for a broad range of objective functions, including those that may be time-varying or switch among a finite set of candidates. %, characterized by mild technical assumptions.

\vspace{0.1cm}
\noindent 
(b) \emph{Safety and Optimality via Hard and Soft Constraints:} The proposed algorithms can satisfy safety-critical constraints \emph{at all times} by using continuous or discontinuous projection maps. The systematic incorporation of these mappings into ES vector fields remained mostly unexplored in the literature, and our results show that they can be safely used in feedback loops to solve optimization problems with hard constraints. In the context of ES, to allow for enough exploration via dithering, the projection maps are applied to a shrunken feasible set that can be made arbitrarily close to the nominal feasible set by decreasing the amplitude of the dithers. In this way, the algorithms are able to provide suitable evolution directions near the boundary of the feasible set, achieving a property of ``practical safety'', similar in spirit to the one studied in \cite{safeES}. In addition to the hard constraints, the proposed controllers are also able to simultaneously handle soft constraints via primal-dual ES vector fields, thus achieving safety and optimality in a variety of model-free optimization problems. 
    
\vspace{0.1cm}
\noindent 
(c) \emph{Stability and Performance Guarantees:} We leverage averaging and singular perturbation theory for non-smooth (and hybrid) systems, as well as Lyapunov-based arguments, to show that the proposed dynamics can guarantee convergence to an arbitrarily small neighborhood of the optimal set, from arbitrarily large compact sets of initial conditions in the feasible set. Moreover, by exploiting the well-posedness of the dynamics and the optimization problem, the algorithms also guarantee suitable robustness properties with respect to small bounded additive disturbances acting on the states and dynamics of the closed-loop system. This is a fundamental property for practical applications and  is non-trivial to achieve in model-free algorithms. We also provide tracking bounds for time-varying optimization problems using (practical) input-to-state stability tools, and we provide stability results for a class of ES problems  with unknown switching objective functions,  which have remained mostly unexplored in the literature.
%
% \vspace{0.1cm}
% \noindent     
% (d) \emph{Decentralized Implementation:} The P-PDZD are extended to cooperative multi-agent problems, leading to decentralized algorithms where each agent 
%     receives broadcast feedback from its neighbors and the unknown environment, and only needs to perform local computations and local actions. As a consequence, there is no need for a centralized center that collects global information on the network, thus preserving the data privacy of each agent.

% \vspace{0.1cm}
% \noindent 
% (d) \emph{Application on Safe Model-free Voltage Control:} We demonstrate the optimality, robustness, and adaptivity of one of the proposed algorithms through extensive numerical experiments on a real-world problem: safe model-free voltage control in a power distribution network. %\lina{I need to check the numerical section to come back to here}

\vspace{0.1cm}
Earlier, partial results of this paper appeared in the conference paper \cite{chen2021safe}. The results of \cite{chen2021safe} are dedicated only to a particular optimal voltage control problem in power systems using only one of the algorithms studied in this paper. In contrast to \cite{chen2021safe}, in this paper, we consider a generic constrained optimization problem and we study two different types of projection maps (continuous and discontinuous), which require different analytical tools and lead to two different algorithms. Additionally, we present novel tracking results for time-varying optimization problems and switching cost functions, and we establish robustness guarantees for all the algorithms. Unlike \cite{chen2021safe}, we also present the complete proofs of the results, as well as novel illustrative examples.

\vspace{0.1cm}
The remainder of this paper is organized as follows: 
Section \ref{sec:preliminary} introduces the notation and the preliminaries. Section \ref{sec:problem} presents the problem formulation. Section \ref{sec:lipschitz} introduces the projected ZO dynamics that incorporate Lipschitz projection maps, and establishes results for static maps, time-varying maps, and switching maps. Section \ref{sec:discontinuous} considers projected gradient-based ZO dynamics with discontinuous projections. Section \ref{section_analysis} presents the analysis and proofs. Numerical experiments are presented throughout the paper to illustrate the main ideas and results. The paper ends with conclusions presented in Section \ref{sec:conclusion}.
%
%

% \section{PRELIMINARIES}\label{sec:preliminary}
\section{NOTATION AND PRELIMINARIES}\label{sec:preliminary}
%
%This section introduces the main notations and definitions that are used in this paper, and presents the preliminaries on extremum seeking (ES) control.

\vspace{0.1cm}
\subsection{Notation}
We use unbolded lower-case letters for scalars and bolded lower-case letters for column vectors. % $\dot{x}:=\frac{d x}{d t}$ denotes the  derivative of $x$ to time $t$. 
We use $\mathbb{R}_+:=[0,+\infty)$ to denote the set of non-negative  real values and use $\B$ to denote a closed unit ball of appropriate dimension.  %$|\cdot|$  denotes the cardinality of a set. 
%$\mathcal{B}(\bx,\epsilon)$ denotes the open ball centered at $\bx\in\R^n$ with radius $\epsilon>0$.
We use
$||\cdot||$ to denote the Euclidean norm of a vector and  use $[\bx; \by] := [\bx^\top, \by^\top]^\top$ to denote the column merge of column vectors $\bx,\by$. Given a positive integer $n$, we   define the index set $[n]:=\{1,\cdots,n\}$. The distance between a point $\bx\in \R^n$ and a nonempty  closed convex set $\sX\subseteq \R^n$ is denoted as $||\bx||_{\sX}:=\underset{\mathbf{\by}\in \sX}{\inf}\, ||\by-\bx||$; and the Euclidean projection of $\bx$ onto the set $\sX$ is defined as 
\begin{equation}\label{projectionmap}
\mathcal{P}_{\sX}(\bx):= \underset{\by \in \sX}{\mathrm{arg\, inf}}  \, ||\by -\bx||.
\end{equation}
The \emph{norm cone} to a set $\sX$ at a point $\bx\in\sX$ is defined as 
\begin{align} \label{eq:normcone}
    N_{\sX}(\bx):=   \left\{ \bs\in \mathbb{R}^n:\,  \bs^\top( \by-\bx) \leq 0, \ \forall~\by\in \sX  \right\}.
\end{align}
The \emph{tangent cone} to    $\sX$ at a point $\bx\in\sX$ is defined as 
\begin{align} \label{eq:tangcone}
    T_{\sX}(\bx):=  \left\{ \bd\in \mathbb{R}^n:\, \bd^\top \bs\leq 0, \,\forall~ \bs\in N_\sX(\bx)  \right\},
\end{align}
which is the polar cone of the normal cone $N_{\sX}(\bx)$. 
%
% \begin{figure}[t!]
%     \centering
% \includegraphics[width=0.2\textwidth]{TangentCones.eps}~~~~\includegraphics[width=0.24\textwidth]{NormalCones.eps}
%       \caption{{\small Illustration of tangent and normal cones of a set $\mathcal{X}$. {\color{red} Since the definitions are standard, we may remove the figure at this time to save space.}}}
%     \label{fig:power}
%     \vspace{-0.3cm}
% \end{figure}
%
A continuous function $\beta(r,s): \R_+\times \R_+\to \R_+$ is said to be of class-$\mathcal{KL}$ if it is zero at zero, non-decreasing in its first argument $r$, non-increasing in the second argument $s$, $\lim_{r\to0^+}\beta(r,s)=0$ for each $s$, and $\lim_{s\to\infty}\beta(r,s)=0$ for each $r$ \cite[Def. 3.38]{Goebel:12}.

In this paper, we consider constrained dynamical systems given by
\begin{equation}\label{ODE00}
\mathbf{x}\in C,~~~\dot{\mathbf{x}}\in F(\mathbf{x}),
\end{equation}
where $\mathbf{x}\in\mathbb{R}^n$ is the state, $C$ is the    {flow set}, and $F:\mathbb{R}^n\rightrightarrows\mathbb{R}^n$ is the  {flow map}, which can be set-valued. We use $\dot{\mathbf{x}}=\frac{d\mathbf{x}(t)}{dt}$ to denote the time derivative of the function $t\mapsto \mathbf{x}(t)$. A function $\bx$ is said to be a (Caratheodory) solution to \eqref{ODE00} if
1) $t\mapsto \mathbf{x}(t)$ is absolutely continuous on each compact sub-interval of its domain $\text{dom}(\mathbf{x})$; 2) $\mathbf{x}(0)\in {C}$; and 3) $\dot{\mathbf{x}}(t)\in F(\mathbf{x}(t))$ and $\mathbf{x}(t)\in C$ for almost all $t\in\text{dom}(\mathbf{x})$ \cite[pp. 4]{filippov2013differential}. The solution $\mathbf{x}$ is said to be complete if $\text{dom}(\mathbf{x})=[0,\infty)$. 
 If the flow map $F$ is single-valued,  \eqref{ODE00} reduces to an ordinary differential equation. If  $F$ is also continuous, solutions $\mathbf{x}$ to \eqref{ODE00} are continuously differentiable functions. In addition, if $F$ is locally Lipschitz, then solutions to \eqref{ODE00} are unique.
% The stability and convergence properties of our algorithms will be studied with respect to compact sets. To study these properties, we will use class $\mathcal{K}\mathcal{L}$ functions. A continuous function $\beta(r,s): \R_+\times \R_+\to \R_+$ is said to be of class-$\mathcal{KL}$ if it is zero at zero, strictly increasing in the first argument $r$, non-increasing in the second argument $s$ and converging to zero as $s\to +\infty$. 
%\begin{definition}
%(Lipschitz Function)\normalfont \cite{cortes2008discontinuous}. A function $\mathbf{f}:\R^n\to\R^m$ is \emph{locally Lipschitz} at $\bx\in\R^n$ if there exist constants $
%L_{\bx},\epsilon\in(0,\infty)$ such that 
%\begin{align*}
%    ||\mathbf{f}(\by_1) -\mathbf{f}(\by_2) || \leq L_{\bx} ||\by_1-\by_2||,
%\end{align*}
%for all $\by_1,\by_2\in \mathcal{B}(\bx,\epsilon)$. A function is \emph{locally Lipschitz on a set} $\sX \subseteq \R^n$ if it is locally Lipschitz at $\bx$ for all $\bx\in\sX$.  A function $\mathbf{f}:\R^n\to\R^m$  is \emph{globally Lipschitz} on a set $\sX \subseteq \R^n$ if there exists a constant $L$ such that
%\begin{align*}
%\quad     ||\mathbf{f}(\bx_1) -\mathbf{f}(\bx_2) || \leq L ||\bx_1-\bx_2||,\quad \forall \bx_1,\bx_2\in\sX.
%\end{align*}
%\end{definition}
%\lina{to make it more rigorous, would like to provide the equation explaining ``locally Lipschitz gradient''}
%

\vspace{-0.1cm}
\subsection{Preliminaries on Extremum Seeking Control}\label{sec:es:pre}
%
%\red{we need to review the introduction of ES or we can put it to the appendix for a broader audience. i) explain where equation (1) is from--why the right hand side has ``gradient'' approximation; exploration, what's $\omega$, etc; ii) explain the true math behind it by mentioning average system, perturbation, etc. } \red{To analyze the dynamics of this system, we can study the \textit{average} of this system and view system \eqref{eq:simes} as a small}
%

Extremum Seeking (ES) control is a type of adaptive control that is able to steer a plant towards a state that optimizes a particular steady-state performance metric using real-time output feedback. These types of controllers can be seen as continuous-time ZO optimization algorithms with (uniform) convergence and stability guarantees. %Figure \ref{fig:es} shows a classic ES block diagram \cite{Krstic2000,tan2006non}, which is the basis for other more advanced architectures \cite{durr2014extremum,PovedaTAC17,ESAntiwindup2,Nesic_ChineseCDC,ghaffari2012multivariable,scheinker2017model}. 
To explain the rationale behind these algorithms, we consider the optimization problem 
\begin{equation}\label{simpleproblem}
\min_{x} f(x),
\end{equation}
where $f:\mathbb{R}\to\mathbb{R}$ is a function that is at least twice continuously differentiable.  
% averaging-based ES controller \cite{Krstic2000,tan2006non}, which is a continuous-time ZO algorithm that seeks to maximize the plant $f$ (here, assumed to be static), using only output measurements. This scheme
% To understand this scheme,
%ES control  is a type of adaptive control  that utilizes only output feedback to
%steer a dynamical system to a state that  attains an extremum \cite{ariyur2003real}. 
%Hence, ES can be viewed as a continuous-time zeroth-order method to solve optimization problems, which essentially estimates  the gradient of the objective function based on exploratory probing signals. consider the problem of solving 
%$\min_{\mathbf{x}} f(\mathbf{x})$, where $f$ is a suitable cost function with bounded level sets in a neighborhood of the set of minimizers. 
A standard approach to finding the minimizer of $f$ is to use a gradient descent flow in the form $\dot{ {x}} = -k_x\cdot \frac{d f(x)}{d x}$, where the gain $k_x$ defines the rate of evolution of the system. However, when the derivative of $f$ is unknown, gradient flows cannot be directly implemented, and instead, \emph{model-free} techniques are required. To address this issue, ES approximates the behavior of the gradient flow by adding a high-frequency periodic probing signal $\varepsilon_a \hat{\mu}(t)$ with amplitude $\varepsilon_a$ to the nominal input of the plant. The resulting output $y = f(x+ \varepsilon_a \hat{\mu}(t))$, which is assumed to be available for measurements, is then multiplied by the same probing signal $\hat{\mu}(t)$, and further normalized by the constant $2/\varepsilon_a$. The loop is closed with an integrator with a negative gain $-k_x$, leading to the ES dynamics: 
\begin{align}\label{eq:simes}
    \dot{x} = -k_x \frac{2}{\varepsilon_a}  f\left(x+\varepsilon_a \hat{\mu}(t)\right) \hat{\mu}(t).
\end{align}
When the frequency of $\hat{\mu}(\cdot)$ is sufficiently large compared to the rate of evolution $k_x$, 
%A key  fact is  that the ES dynamics \eqref{eq:simes} 
%with  small  $a$ and large $\omega$  behaves approximately like the gradient descent flow  $\dot{x}=-k\cdot\nabla f(x)$, which can steer  $x$ to a (local) minimum $x^*=\arg\min_{x}f(x)$ under appropriate conditions on $f(\cdot)$. Note that to implement the ES dynamics \eqref{eq:simes}, one does not need the knowledge of function $f$ but only its values, and thus ES control is regarded as a zeroth-order method.
the ES dynamics \eqref{eq:simes} exhibits a time scale separation property that allows to approximate the behavior of $x$ based on the average of the vector field of \eqref{eq:simes}. For example, consider the use of a sinusoidal signal as the probing single, i.e.,  $\hat{\mu}(t) := \sin (\omega t)$. With large $\omega>0$ and small $\varepsilon_a$,
%
%, i.e., letting 
%\begin{align*}
%    d(\omega t) := \sin (\omega t), \text{ and } \eta_d=2.
%\end{align*}
% The fast-timescale variation of the ES dynamics   is caused by the periodic sinusoidal signal $\sin(\omega t)$, while the slow-timescale variation  governed by the gain $k$  dominates the evolution of $x$. 
%By averaging theory, the fast variation can be washed out and
%one can obtain a time-invariant average dynamics that describes the main  trend of the evolution  of $x$. Specifically,
%
we  consider the Taylor expansion of $f$: $$ f(x+\varepsilon_a\sin(\omega t)) =  f(x)
+ \varepsilon_a\sin(\omega t)\frac{d f(x)}{d x} +   \mathcal{O}(\varepsilon_a^2).$$ By  computing the average of the vector field of \eqref{eq:simes} over one period $T=\frac{2\pi}{\omega}$ of the probing signal, one obtains
%\begin{align}
%\begin{split}
%    f(x+a&\sin(\omega t)) \sin(\omega t) \\
%=& f(x) \sin(\omega t)
%+ a\sin(\omega t)^2\frac{\partial f(x)}{\partial x} +   %\mathcal{O}(a^2).
%\end{split}
%\end{align}
%
%Thus the average dynamics of (\ref{eq:simes}) is given by 
%
\begin{align}\label{eq:avedyn}
 \dot{{x}}& = \frac{1}{T} \int_{0}^T\! \!-k_x\frac{2}{\varepsilon_a}  f(x+\varepsilon_a\sin(\omega t))\sin(\omega t)\,  dt \nonumber\\
 &
= -\frac{k_x}{T} \int_{0}^T\! \!2\sin^2(\omega t)\frac{d f(x)}{d x} + \sO(\varepsilon_a) \, dt , \nonumber\\
& = -k_x\frac{d f(x)}{d x} +  \sO(\varepsilon_a):=h_\text{ave}(x)
 \end{align}
 where $\sO(\varepsilon_a)$ denotes high-order terms, bounded on compact sets, that vanish as $\varepsilon_a\to0^+$. The average system \eqref{eq:avedyn} is essentially an $\sO(\varepsilon_a)$-perturbed gradient descent flow. Under suitable assumptions on $f$, averaging theory and perturbation theory show that the trajectories of \eqref{eq:simes} will approximate those of \eqref{eq:avedyn} (on compact sets and compact time intervals) as $\varepsilon_a\to0^+$ and as $\omega\to0^+$ \cite[Theorem 1]{teel2003unified}. Uniform stability properties of gradient flows can then be leveraged to establish stability results for \eqref{eq:simes} in the infinite horizon \cite[Theorem 2]{teel2003unified}. This analysis can also be applied to the multi-variable case using an appropriate choice of the (vector) frequencies $\omega$, and to other architectures using Lie-bracket averaging theory that results in similar average systems \cite{DurrLieBracket,scheinker2017model}. %The majority of the ES results and architectures consider only unconstrained optimization problems, optimization problems with asymptotic constraints, or problems defined on Euclidean manifolds.  

\section{PROBLEM FORMULATION} \label{sec:problem}
%
%In this section, we define and motivate our main problem of interest, namely, model-free optimization with hard and soft constraints; as well as the different settings defined by the available information.
%
%the setting of available information. Then, we present several  practical application examples to explain the motivation and justify the problem setting. After that, we solve the problem  using projected primal dual gradient dynamics (P-PDGD) methods. In particular, two types of P-PDGDs with continuous and discontinuous projections are considered.
%
%\vspace{-0.3cm}
%\subsection{Problem Formulation} \label{sec:problem:intro}
%
In contrast to \eqref{simpleproblem}, in this paper, we consider \emph{constrained} optimization problems of the form
\begin{subequations}\label{eq:gen}
\begin{align}
      \text{Obj.} \ \, &  \min_{\bx}\, f(\bx) \label{eq:gen:obj}\\
    \text{s.t.} \ \    &\ \bx\in \sX \label{eq:gen:set}\\
      &\ g_j(\bx)\leq 0, \qquad j\in[m],\label{eq:gen:ineq}
 \end{align}
\end{subequations}
where $\bx \in\R^n$ is the decision variable, $f:\R^n\to\R$ is the  objective function, $\sX\subseteq \R^n$ denotes the feasible set of $\bx$, and the vector-valued function $\bg:=[g_1;g_2;\cdots;g_m]:\R^n\to\R^m$ describes additional inequality constraints on $\bx$. The set of the optimal solutions of \eqref{eq:gen} is denoted as $\mathcal{X}^*\subset\mathbb{R}^n$.

  \textit{Information Availability:}
We consider  the problem setting where the feasible set $\mathcal{X}$ is known but the mathematical forms of $f(\cdot)$ and $\bg(\cdot)$ are unknown. In this case, one can only query (in real-time) the  values of $f(\bx)$ and $\bg(\bx)$ for a given $\bx$. That is, the optimization solver can only access the zeroth-order information of $f(\cdot)$ and $\bg(\cdot)$, but not their (first-order) gradients  or (second-order) Hessian information.

%\begin{remark}\normalfont
The motivation and rationale of the above problem setting are explained below:

\begin{itemize}  
    \item [1)] The above problem is motivated by the  feedback control design that seeks to steer an unknown plant in real time to an optimal solution of problem \eqref{eq:gen}. Here, we model the plant using the static input-to-output maps $f(\cdot)$ and $\bg(\cdot)$ to approximate its steady-state response. The validity of this approximation lies in the fact that in many applications the plant is a stable dynamical system that converges to a steady state in a much faster time scale compared to the controller. The steady-state approximation of the plant can then be justified using a singular perturbation argument \cite[Theorem 2]{teel2003unified}, provided the time-scale separation is sufficiently large. %guarantee that the steady-state approximation of the plant is valid provided the time scale separation in the closed-loop system is sufficiently large.
    %the intrinsic time scale separation between the dynamics of the plant and the controller, 
    %
    %which are assumed to be sufficiently fast so that they
    %
    %is as we consider fast stable plant dynamics that converge immediately given any input $\bx\in\sX$ and aim  to optimize the steady-state performance.
    
     %It is   referred to as the timescale separation property 
     % with fast plant dynamics and  slow control, which 
       %is commonly assumed in ES control \cite{tan2010extremum,ariyur2003real}.
 %   The complete dynamical model of the plant can be formulated as 
%\begin{subequations}
%\begin{align}
 %  \dot{ \bz} & = \bh(\bz,\bx) \label{eq:sim:state}\\
  %  \by & = [\bar{f}(\bz,\bx)); \bar{\bg}(\bz,\bx)],\label{eq:sim:output}
%\end{align}
%\end{subequations}
 %   where $\bz$ is the system state, $\bx$ is the input, and $\by$ is the output.    
%
\item [2)] For many complex engineering systems, their  models, captured by the maps $f(\cdot)$ and $\bg(\cdot)$, may be unknown, unavailable, or too costly to estimate. 
   % optimal control of various cyber-physical systems with the goal of steering the decision $\bx$ to an optimal solution of  problem      \eqref{eq:gen}, while the system models (described by functions $f$ and $\bg$) are unavailable or too costly to procure in practice. 
On the other hand, the widespread deployment of smart meters and sensors provides real-time  measurements of the system outputs. These measurements can be 
 interpreted as the function evaluations of $f(\cdot)$ and $\bg(\cdot)$ and  can be used as the system feedback to circumvent the unknown model information. 
\item [3)] In problem \eqref{eq:gen}, we
    distinguish  \emph{hard constraints}, modeled by $\sX$, and \emph{asymptotic constraints}, modeled by the inequalities  \eqref{eq:gen:ineq}. Thus, the constraints imposed by $\sX$  \eqref{eq:gen:set} must be satisfied at all times, while inequalities \eqref{eq:gen:ineq} 
    may be violated during the transient process but should be  satisfied   in the steady states.
\end{itemize}

This paper aims to develop model-free feedback optimization algorithms that are able to solve problem \eqref{eq:gen} using only zeroth-order information, while simultaneously  satisfying hard and asymptotic constraints. To achieve these goals, in Sections \ref{sec:lipschitz} and \ref{sec:discontinuous}, we will study a class of ZO feedback optimization algorithms that are based on ES and incorporate two types of projection maps. To guarantee that problem \eqref{eq:gen} is well-posed, throughout this paper we will make the following assumptions, which, as discussed later, can be used to relax standard global convexity assumptions considered in the literature of ES.
\vspace{1pt}
\begin{assumption} \label{ass:con_sm}
 The feasible set $\mathcal{X}$ is  nonempty, closed, and convex. The functions $f$ and $g_1,\cdots,g_m$ are convex and at least twice continuously differentiable on an open set containing $\mathcal{X}$. The function $f$ is radially unbounded.%, whose gradients $\nabla_{\bx} f$ and $\nabla_{\bx}g_1,\cdots,\nabla_{\bx}g_m$  are  locally Lipschitz  on $\mathcal{X}$
\QEDB
\end{assumption}
%

%\vspace{0.1cm}
% To guarantee the well-posedness of problem \eqref{eq:gen}, we also impose the following assumption.
%

%\vspace{0.1cm}
\vspace{1pt}
\begin{assumption} \label{ass:finite}
Problem (\ref{eq:gen}) has a finite optimum and the Slater's conditions hold. Moreover, the set  of optimal solutions $\sX^*$ is compact. \QEDB
\end{assumption}
%

%\vspace{0.1cm}
%Assumptions \ref{ass:con_sm} and \ref{ass:finite} are mainly used for  theoretical analysis, but our proposed algorithms can be practically applied to a wider range of problems that may not satisfy these assumptions. % globally, but only in a region of the state space, c.f. Figure \ref{fig:regionalresult}. This observation will be further validated by the numerical experiments in Section \ref{sec:simulation}.

% \input{Application}

\section{MODEL-FREE FEEDBACK OPTIMIZATION  WITH LIPSCHITZ PROJECTIONS}
\label{sec:lipschitz}
In this section, we introduce a class of gradient-based continuous-time ZO algorithms that incorporate Lipschitz continuous projection maps. We term these algorithms as the \emph{projected gradient-based zeroth-order} (P-GZO) dynamics. We first study a reduced version of \eqref{eq:gen} that considers only the hard constraint \eqref{eq:gen:set}, i.e., we consider the problem:
\begin{align}\label{eq:gen:reduce}
       \min_{\bx\in \sX} \ f(\bx).
 \end{align}
For this problem, we establish stability, safety, and  tracking results under the P-GZO dynamics. After this, we develop results for the case when $f$ is dynamically drawn from a finite collection of cost functions that share the same minimizer, a problem that emerges in systems with switching plants or costs. Lastly, we further incorporate the constraints \eqref{eq:gen:ineq} using a projected primal-dual zeroth-order (P-PDZO) algorithm. 
%To address problem \eqref{eq:gen}, we first consider a class of zeroth-order (ZO) algorithms that incorporate  projection maps that are Lipschitz continuous. We first focus on establishing stability and safety results for static model-free optimization problems with constraints of the form \eqref{eq:gen:set}. After this, and under stronger assumptions on the cost function, we establish tracking results via ISS bounds for the case when the (unique) solution of \eqref{eq:gen} changes over time. Subsequently, we incorporate the soft constraints \eqref{eq:gen:ineq} using a model-free primal-dual formulation.
%
% with the projected primal-dual gradient dynamics (P-PDGD) method which uses the gradient information of $f$ and $g$.
% Then, we integrate ES control into P-PDGD and 
% develop the projected primal-dual zeroth-order dynamics (P-PDZD) method, which is the model-free feedback solution algorithm for solving \eqref{eq:gen}. Lastly, we propose the decentralized application of P-PDZD for cooperative multi-agent systems.
%

%\vspace{-0.2cm}
%

\subsection{GZO Dynamics with Lipschitz Projection} 
\label{sec:pzodesign}
To solve \eqref{eq:gen:reduce}, we consider the following dynamics, termed the
\emph{projected gradient zeroth-order} (P-GZO) dynamics:
\begin{subequations}\label{eq:pzo}
\begin{align}
\dot{\mathbf{x}}&=k_x \Big(\mathcal{P}_{\mathcal{X}}(\mathbf{x}-\alpha_x \,\bm{\xi})-\mathbf{x}\Big),\label{eq:pzo:x}\\
\dot{ {\bxi}}&=\frac{1}{\varepsilon_{\xi}}\Big(-{\bxi}+\frac{2}{\varepsilon_a}f(\hat{\mathbf{x}}) \hat{\bmu}\Big),\label{eq:pzo:xi} \\
\dot{ {\bmu}} &=\frac{1}{\varepsilon_{\omega}}\Lambda_{\kappa} {\bmu},\label{eq:pzo:mu}
\end{align}
\end{subequations}
where $k_x,\alpha_x,\varepsilon_{\xi},\varepsilon_a,\varepsilon_\omega\!>\!0$ are tunable parameters. The dynamics \eqref{eq:pzo:x} incorporates a Lipschitz projection map of the form \eqref{projectionmap} to ensure that $\bx$ stays within the feasible set $\sX$. The dynamics \eqref{eq:pzo:xi}  estimates the gradient $\nabla f$ with a new state $\bxi\in\R^n$, whose dynamics depend on the \emph{measured} output $y=f(\hat{\bx})$, where $\hat{\bx}$ is the perturbed input   defined as
\begin{align} \label{eq:hatx}
    \hat{\bx} := \bx  + \varepsilon_a \hat{\bmu}.
\end{align}
In \eqref{eq:hatx}, $\hat{\bmu}:\mathbb{R}_{\geq0}\to \R^n$ is a vector-valued  periodic dither signal that is  generated by the linear dynamic oscillator \eqref{eq:pzo:mu}. Specifically, the vector $\hat{\bmu}$ collects all the odd entries of the state $ {\bmu}\in \R^{2n}$, i.e.,
\begin{align} \label{eq:hatmu}
    \hat{\bmu} & := [ {\mu}_1,  {\mu}_3,  {\mu}_5,\cdots, {\mu}_{2n-1}]^\top.
\end{align}
The matrix $\Lambda_{\kappa}\in\mathbb{R}^{2n\times 2n}$ in \eqref{eq:pzo:mu} is block diagonal, with the $i$-th diagonal block given by
\begin{equation}\label{eq:lambdaimaatrix}
\Lambda_{\kappa_i}=\left[\begin{array}{cc}
0 & -2\pi\kappa_i\\
2\pi\kappa_i & 0
\end{array}\right]\in \R^{2\times 2}, \quad i\in[n],
\end{equation}
which is parameterized by the tunable constant $\kappa_i>0$. Hence, \eqref{eq:pzo:mu} describes $n$ autonomous oscillators, whose solutions $ {\bmu}$ can be explicitly computed as
 \begin{align} \label{eq:periodicdither}
\begin{split}
       {\mu}_i(t) = {\mu}_{i}(0)\sin\left(\frac{2\pi\kappa_i}{\varepsilon_{\omega}} t\right)&
+\mu_{i+1}(0)\cos\big(\frac{2\pi\kappa_i}{\varepsilon_{\omega}} t\big),   \\
&   \qquad \forall\, i = 1, 3, \cdots, 2n\!-\!1,
  \end{split}
\end{align}  
and we choose initial conditions that lie on the unit circle: 
\begin{equation}\label{initialization}
{\mu}_{i}(0)^2 + {\mu}_{i+1}(0)^2=1. 
\end{equation}
%
%By  definition \eqref{eq:hatmu}, the   dither signal $\hat{\bmu} (t)$ is sinusoidal and periodic, whose entries are given by  \eqref{eq:periodicdither}. 
For example, when $ {\mu}_i(0) = 1$ and $ {\mu}_{i+1}(0) = 0$ for $i = 1, 3, \cdots, 2n\!-\!1$, equation \eqref{eq:hatmu} becomes
$$\hat{\bmu}(t):= \Big[\sin\left(\frac{2\pi\kappa_1}{\varepsilon_{\omega}} t\right), \cdots, \sin\left(\frac{2\pi\kappa_{2n\!-\!1}}{\varepsilon_{\omega}} t\right) \Big]^\top.$$

In addition to sinusoidal dither signals, other types of dither signals can also be employed to obtain suitable estimations of the gradient, including triangular waves and square waves, see \cite{tan2008choice,scheinker2016bounded,poveda2018hybrid}.  By incorporating the linear dynamic oscillator \eqref{eq:pzo:mu}, the P-GZO dynamics \eqref{eq:pzo} becomes an autonomous system, which facilitates the theoretical analysis. 

\begin{figure}[t!]
    \centering
\includegraphics[width=0.45\textwidth]{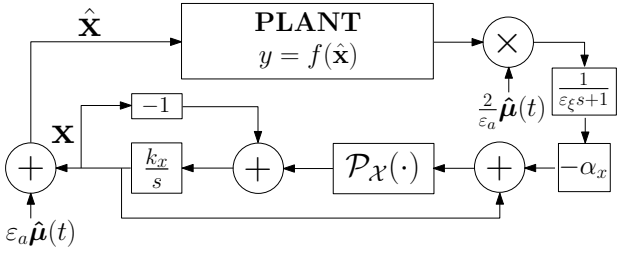}
      \caption{\small{Block diagram of P-GZO dynamics.}}
    \label{fig:Fig0}
    \vspace{-0.3cm}
\end{figure}

The P-GZO dynamics \eqref{eq:pzo}, with the overall state $\bz:=(\bx, \bxi,  {\bmu})$, are defined with respect to the following flow set
\begin{equation}\label{flowset1}
\mathbf{C}_1:=\mathcal{X}\times \mathbb{R}^n\times\mathbb{T}^n,
\end{equation}
where $\mathbb{T}^n:=\mathbb{S}\times\mathbb{S}\times\ldots\times\mathbb{S}$ %is the $n^{th}$ torus. 
and $\mathbb{S} \subset \R^2$ denotes the unit circle centered at the origin. By construction and Assumption \ref{ass:con_sm}, the set $\mathbf{C}_1$ is closed, and it enforces condition \eqref{initialization} on the initialization of the state $\bmu$. Note that the P-GZO dynamics \eqref{eq:pzo} has a Lipschitz continuous vector field on the right-hand side due to the use of a Lipschitz projection mapping. Figure \ref{fig:Fig0} shows a block diagram of the proposed algorithm.

%We term the constrained ODE \eqref{eq:pzo}-\eqref{flowset1} as the \emph{projected gradient zeroth-order} (P-GZO) dynamics, and we note that this system has a Lipschitz continuous vector field on the right-hand side due to the use of a Lipschitz projection mapping. Figure \ref{fig:Fig0} shows a block diagram of the proposed algorithm.

The following assumption will be used throughout  this paper to distinguish different dither signal components with different frequency parameters.
\vspace{1pt}
\begin{assumption}\label{ass:frequenciesdithers}
The parameters $\kappa_i>0$ in \eqref{eq:lambdaimaatrix} are rational numbers and satisfy $\kappa_i\neq \kappa_j$, and $\kappa_i\neq 2\kappa_j$ for all $i\neq j$. \QEDB
\end{assumption}

\vspace{0.1cm}
%
%
%\begin{subequations}\label{periodicdither}
 % \begin{align}
%\hat{\mu}_i(t) &=\hat{\mu}_{i}(0)\sin\big(\frac{2\pi\kappa_i}{\varepsilon_{\omega}} t\big)\!+\!\mu_{i\!+\!1}(0)\cos\big(\frac{2\pi\kappa_i}{\varepsilon_{\omega}} t\big), \\
%\hat{\mu}_{i\!+\!1}(t) &= - \hat{\mu}_{i}(0)\cos\big(\frac{2\pi\kappa_i}{\varepsilon_{\omega}} t\big)\!+\!\mu_{i\!+\!1}(0)\sin\big(\frac{2\pi\kappa_i}{\varepsilon_{\omega}} t\big),
%\end{align}  
%\end{subequations}
%for $i\in\{1, 3, \cdots, 2n\!-\!1\}$. 

%The first continuous-time zeroth-order algorithm that we consider makes use of an auxiliary input state $\mathbf{x}\in\mathbb{R}^n$, whose dynamics incorporate a Lipschitz projection map, a state $\mathbf{\xi}$, generated by a low-pass filter, and a dynamic oscillator with state $\mathbf{\mu}$. The complete equations describing the algorithm are
%

We further explain the proposed P-GZO dynamics \eqref{eq:pzo} with the following remarks.

\vspace{2pt}
\begin{remark}\label{rem:lowpass}
 \tcb{The intuition behind the P-GZO dynamics \eqref{eq:pzo} is that, for sufficiently small values of $\varepsilon_{\omega}$ and $\varepsilon_{a}$, the term $\frac{2}{\varepsilon_a}f(\hat{\mathbf{x}})\hat{\bmu}$ provides, on average, an $\sO(\varepsilon_a)$-approximation of the gradient $\nabla f(\bx)$. Similar one-point estimation mechanisms are common in the literature of zeroth-order methods and stochastic approximations via simultanous perturbations \cite{stochasticESC}, although in our case the dithers are deterministic. The dynamics \eqref{eq:pzo:xi} with a small $\varepsilon_\xi$ behaves as a low-pass filter with input  $\frac{2}{\varepsilon_a}f(\hat{\mathbf{x}})\hat{\bmu}$ and output $\bxi$, which, at steady state, satisfies $\bxi=\nabla f(\bx)+\mathcal{O}(\varepsilon_a)$. This filter is tuned to operate in a faster time scale compared to the dynamics of $\mathbf{x}$. In this way, the low-pass filter facilitates the analysis of the projected system via averaging theory by removing from $\dot{\bx}$ the term that explicitly includes the highly oscillatory signal $\hat{\bmu}$. Otherwise, the projection in \eqref{eq:pzo:x} may interfere with the computation of the average dynamics of $\bx$ near the boundary of $\sX$.} %See the analysis in Section {\color{red}add section NO.} for more explanations.
%the signal $f(\hat{\mathbf{x}})\mathbf{d}(\mu)$, used directly in \eqref{eq:simes} as a surrogate of the gradient of $f$, is now the input of a low-pass filter with state $\xi$. As shown in the analysis, presented in Section \ref{section_analysis}, this extension is needed to guarantee that, on average, $\dot{\mathbf{x}}$ evolves approximately following a projected gradient flow under suitable tuning of the parameters of the system. Indeed, if the filter is not used, the projection in \eqref{projectionX} might interfere with the computation of the average near the boundary of $\mathcal{X}$. \QEDB
\end{remark}

\vspace{2pt}
\begin{remark} \label{rem:safety}
%The intuition behind the dynamics \eqref{algo1} is the following: as $\varepsilon_{\omega}\to0^+$, the high-frequency dithers $\mathbf{d}(\mu)$ will make $f(\hat{\mathbf{x}})\mathbf{d}(\mu)$ a suitable $\mathcal{O}(a)$ approximation of $\nabla f$. 
(\emph{Safety and Optimality}). As we will show below in Lemma \ref{lem:forward}, the projection map in \eqref{eq:pzo:x} guarantees that $\bx$ remains always in the feasible set $\sX$,
and thus the actual decision input $\hat{\bx}$ in \eqref{eq:hatx} remains in a small tunable $\sO(\varepsilon_a)$-neighborhood of $\sX$. This property defines a notion of ``practical" safety, similar to those studied in \cite{safeES,Poveda:15}. However, in contrast to other constrained model-free algorithms that use barrier functions \cite{Poveda:20TAC}, orthogonal projections \cite{Poveda:15}, or safety filters \cite{safeES}, the state $\bx$ in \eqref{eq:pzo:x} can actually converge to the boundary of $\mathcal{X}$ in a finite time, a situation that emerges in problems with saturation constraints. On the other hand, for applications where the decision input $\hat{\bx}$ must stay exactly within $\sX$ for all time, the projection map can be applied to a shrunk feasible set $\mathcal{X}_{\varepsilon_a}$ satisfying $\mathcal{X}_{\varepsilon_a}\!+\!\varepsilon_a\mathbb{B}\subseteq\mathcal{X}$. In addition, as stated below in Theorem \ref{thm:pzo}, when $\varepsilon_\xi$ and $\varepsilon_a$ are also sufficiently small, the trajectory $\bx$ of \eqref{eq:pzo} will converge to a small tunable neighbor of the optimal set $\sX^*$ that solves problem \eqref{eq:gen:reduce}.   
\end{remark}
%
%\vspace{2pt}

%The main advantage of system \eqref{algo1} in comparison to the traditional continuous-time ES dynamics \cite{Krstic2000,tan2006non}, is that the projection operator $\mathcal{P}_{\mathcal{X}}$ will guarantee that the trajectories $t\mapsto\mathbf{x}(t)$ remain in the feasible set $\mathcal{X}$ for all times in the domain of the solution. In cases where $\mathcal{X}$ is bounded, the use of the projection might also simplify the analysis since in this case the solutions would be (uniformly) bounded.  Moreover, since \eqref{algo1} is locally Lipschitz, the existence of solutions is always guaranteed. 

\subsection{Stability Analysis of the P-GZO Dynamics}

%Since the forward invariance of $\mathcal{X}$ will play an important role in this paper, we formalize this property in the following lemma.  
%
% This property is formalized in the following Lemma, adapted from \cite[Thm. 3.2]]{Gao}, for closed and convex sets $\mathcal{X}$, the Lipschitz projection $\mathcal{P}_{\mathcal{X}}(\cdot)$ guarantees that the inputs $\hat{x}$ generated via \eqref{nominal_input} remain in a small inflation of $\mathcal{X}$, whenever they exist. Since \eqref{algo1} is locally Lipschitz, existence of solutions is always guaranteed. 
%

To study the P-GZO dynamics \eqref{eq:pzo}, we first establish the following lemma, which shows that the solutions $\bz$ to the ODE \eqref{eq:pzo} (with flow set $\mathbb{R}^n\times\mathbb{R}^n\times\mathbb{T}^n$) remain in $\mathbf{C}_1$ for all time. The proof is presented in Appendix \ref{app:lem:forward:pf}. 
 
\vspace{2pt}
\begin{lemma}\label{lem:forward} %\label{lemma1}
Suppose that Assumption \ref{ass:con_sm} holds.
 Let $\mathbf{z}:=(\mathbf{x},\bm{\xi},\bm{\mu})$ be a solution to \eqref{eq:pzo} with $\mathbf{z}(0)\in \mathbf{C}_1$.   Then,  
 $\mathbf{z}(t)\in \mathbf{C}_1$ and $\hat{\mathbf{x}}(t)\in\mathcal{X}+\varepsilon_a\mathbb{B}$ for all $t\in\text{dom}(\mathbf{z})$. \QEDB
\end{lemma}
\vspace{2pt}

We analyze the stability and convergence properties of the  P-GZO dynamics \eqref{eq:pzo} based on the properties of a nominal ``\emph{target system}'', given by
\begin{equation} \label{eq:nominal:pzo} %\label{projectedgradientflow1}
\mathbf{p}\in \sX,\quad \dot{\mathbf{p}}=k_x\Big(\mathcal{P}_{\mathcal{X}}\big(\mathbf{p}-\alpha_x\nabla f(\mathbf{p})\big)-\mathbf{p}\Big),
\end{equation}
which has been well studied in  the literature \cite{nagurney2012projected}. The following theorem, which is the first result of this paper, only relies on assuming the well-posedness of \eqref{eq:gen:reduce} and suitable stability properties for \eqref{eq:nominal:pzo}. Particular cases where these assumptions are satisfied are discussed afterwards.

\vspace{0.1cm}
\begin{theorem} \label{thm:pzo}%\label{theorem1}
%Consider the constrained ODE \eqref{algo1}-\eqref{flowset1}, and 
Suppose that Assumptions \ref{ass:con_sm}-\ref{ass:frequenciesdithers} hold, and 
\begin{enumerate}[(a)]
\item Every solution of \eqref{eq:nominal:pzo} with $\mathbf{p}(0)\in\mathcal{X}$ is complete;
\item System \eqref{eq:nominal:pzo} renders the optimal set $\mathcal{X}^*$ forward invariant and uniformly attractive. 
\end{enumerate}
Then, for any $\Delta>\nu>0$ there exists $\hat{\varepsilon}_{\xi}>0$ such that for all $\varepsilon_{\xi}\in(0,\hat{\varepsilon}_{\xi})$ there exists $\hat{\varepsilon}_a>0$ such that for all $\varepsilon_a\in(0,\hat{\varepsilon}_a)$, there exists $\hat{\varepsilon}_{\omega}>0$ such that for all $\varepsilon_{\omega}\in(0,\hat{\varepsilon}_{\omega})$,  every solution $\mathbf{z}$ of the P-GZO dynamics 
with $\mathbf{z}(0)\in \mathbf{C}_1\cap \left(\left(\mathcal{W}_1^*+\Delta\mathbb{B}\right)\times\mathbb{T}^n\right)$ is complete and satisfies:
\begin{align}
&\text{({Practical Convergence}):}\quad \limsup\limits_{t\rightarrow\infty}  ||\mathbf{x}(t)||_{\mathcal{X}^*}\leq \nu,\qquad \label{KLbound1}\\
&\text{({Practical Safety}):}\ \ \mathbf{x}(t)\!\in\!\mathcal{X},\ \hat{\mathbf{x}}(t)\!\in\!\mathcal{X}\!+\!\varepsilon_a\mathbb{B},\ \forall t\geq 0, \label{eq:safetybound}
\end{align}
where $\mathcal{W}_1^*:=\{(\mathbf{x},\bxi)\in\mathbb{R}^{2n}:\mathbf{x}\in\mathcal{X}^*,\,\bxi=\nabla f(\mathbf{x})\}$. \QEDB
\end{theorem}

\vspace{0.2cm}
 The complete proof of Theorem \ref{thm:pzo} is presented in Section \ref{sec:pzo:proof} as a particular case of a more general result presented later in Theorem \ref{thm:p-pdzo}. The result of Theorem \ref{thm:pzo} establishes two main properties: 1) convergence from arbitrarily large pre-defined $\Delta$-compact sets of initial conditions to arbitrarily small $\nu$-neighborhoods of the optimal set, which is a typical property of zeroth-order algorithms; and 2) the safety result \eqref{eq:safetybound} for $\bx$ and $\hat{\bx}$ that holds for all time $t\geq 0$. Note that our results do not assume that the feasible set $\mathcal{X}$ is bounded, but, when this is the case, the result becomes global with respect to $\mathcal{X}$.

\begin{figure}[t!]
    \centering
\includegraphics[width=0.48\textwidth]{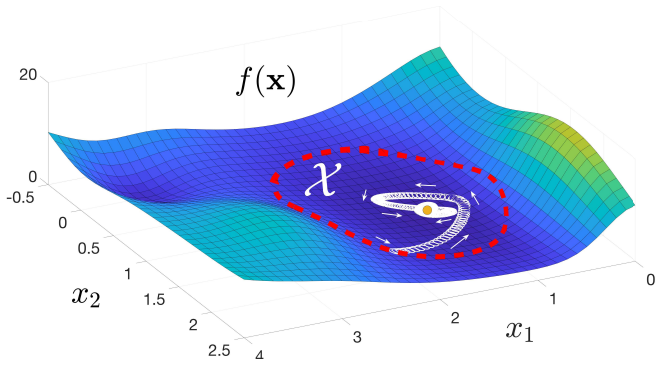}
      \caption{\small{Trajectory $\hat{\mathbf{x}}$ of the P-GZO algorithm on a regionally convex landscape. The safe region $\mathcal{X}$ is delimited by the red dashed line. All trajectories remain in $\mathcal{X}$ and converge to a neighborhood of $\mathbf{x}^*$.}}
    \label{fig:regionalresult}
    \vspace{-0.4cm}
\end{figure}

%Particular cases where these assumptions are satisfied are discussed afterwards. The complete proof is presented in Appendix \ref{sectionproofs} as a particular case of a more general result presented later in Theorem 3. The result of Theorem \ref{theorem1} establishes two main properties: 1) a convergence result that holds \emph{in the limit}, from arbitrarily large $\Delta$-compact sets of initial conditions, to arbitrarily small $\nu$-neighborhoods of the optimal set, which is typical in zeroth-order algorithms; and 2) a safety result for $x$ and $\hat{x}$ that holds \emph{for all times} $t\geq0$. The result does not assume that $\mathcal{X}$ is bounded, but, when this is the case, the semi-global results become global.
%
The conditions under which the assumptions (a) and (b) in Theorem \ref{thm:pzo}  hold for the nominal system \eqref{eq:nominal:pzo} have been extensively studied in the literature \cite{gao2003exponential}. For example, these two assumptions hold when the objective function $f$ is strictly convex \cite[Theorem 1]{gao2003exponential}, in which case $\mathcal{W}_1^*$ is a singleton.  In fact, the result of Theorem \ref{thm:pzo} holds even when $\nabla f$ in \eqref{eq:nominal:pzo} is replaced by a general strictly monotone mapping, since in this case assumptions (a) and (b) of Theorem \ref{thm:pzo} also hold \cite[Corrollary 1]{gao2003exponential}. This implies that Theorem \ref{thm:pzo} can also be used for decision-making problems in games using the pseudo-gradient instead of the gradient, similar to the studies presented in \cite{Frihauf12a}.
\begin{remark}
One of the main limitations of traditional zeroth-order algorithms that emulate gradient descent, such as \eqref{eq:simes}, is that the cost $f$ might not be convex (or gradient-dominated) in the whole space $\mathbb{R}^n$, precluding semi-global convergence results. In this case, projection maps can be used to restrict the evolution of the algorithm to ``safe'' regions $\mathcal{X}$ where suitable convexity/monotonicity properties are presumed to be satisfied. This observation is illustrated in Figure \ref{fig:regionalresult}, where a non-convex landscape, with multiple local minima, maxima, and saddle points, is ``safely'' optimized in a set $\mathcal{X}$ where the assumptions of Theorem \ref{thm:pzo} hold. \QEDB
\end{remark}
\vspace{-0.2cm}
\subsection{Tracking Properties of P-GZO Dynamics}
For many practical applications, the corresponding optimization problem \eqref{eq:gen} is not static but time-varying, with objectives and constraints that may change over time. This subsection considers the time-varying optimization setting by allowing the cost  \( f \) in \eqref{eq:gen:reduce} to depend on a time-varying parameter \( \bth \in \mathbb{R}^p \), i.e., we now consider continuous differentiable mappings $(\bx,\bth)\mapsto f(\bx,\bth)$. In addition,  $\bth$ is assumed to be generated by an (unknown) exosystem of the form
%In many applications, the main advantage of using projection maps in model-free optimization algorithms is their ability to maintain for all times the trajectories $\mathbf{x}$ in the feasible set $\mathcal{X}$, despite sudden and persistent changes in the optimal solutions to problem \eqref{eq:gen}. This time-varying optimization setting can be studied for the model-free dynamics \eqref{algo1} by letting the cost $f$ to depend on a time-varying parameter $\theta\in\mathbb{R}^p$ assumed to be generated by an (unknown) exosystem of the form
%
\begin{equation}\label{eq:exosystem}
\bth\in \Theta,\quad \dot{\bth}\in \varepsilon_{\theta}\Pi(\bth),
\end{equation}
where $\varepsilon_{\theta}>0$ is a parameter that describes the rate of change of $\theta$, $\Theta\subset\mathbb{R}^p$ is a compact set, and $\Pi:\mathbb{R}^p\rightrightarrows\mathbb{R}^p$ is a set-valued mapping assumed to be outer-semicontinuous, locally bounded, and convex valued \cite{Goebel:12}. Additionally, system \eqref{eq:exosystem} is assumed to render the set $\Theta$ strongly forward invariant. By considering exosystems of the form \eqref{eq:exosystem}, we can model a broad family of locally absolutely continuous functions $t\mapsto\theta(t)$. For the case when $\Pi$ is a single-valued mapping, our assumptions are satisfied when $\Pi(\cdot)$ is continuous. As a result, the optimizer $\bx^*$ is also time-varying and describes a trajectory $t\mapsto \mathbf{x}^*(\bth(t))$. In this case, we examine the tracking performance of the P-GZO dynamics in solving the time-varying version of problem \eqref{eq:gen}. We make the following regularity assumptions on the parameterized optimal trajectory.
%In this case, we study the tracking properties of \eqref{algo1} with respect to the time-varying optimizer $t\mapsto \mathbf{x}^*(\theta(t))$, which is guaranteed to be unique under the following additional assumptions:
%
%

\vspace{0.1cm}
\begin{assumption}\label{ass:track}
There exists a continuously differentiable function $\bd:\mathbb{R}^p\to\mathbb{R}^n$ such that
\begin{equation}\label{eq:dtheta}
\mathbf{x}^*(\bth):=\bd(\bth)=\text{arg}\min_{\mathbf{x}\in\mathbf{\mathcal{X}}}f(\mathbf{x},\bth),
\end{equation}
for all $\bth\in\Theta$. Also, there exist $\ell,\gamma>0$ such that 
\begin{subequations}
\begin{align}
&||\nabla f(\mathbf{x},\bth)-\nabla f(\mathbf{y},\bth)||\leq \ell||\mathbf{x}-\mathbf{y}||,\label{eq:globallipschitz1}\\
&f(\mathbf{x},\bth)\!-\!f(\mathbf{y},\bth)\geq \nabla_{\by} f(\mathbf{y},\bth)(\mathbf{x}\!-\!\mathbf{y})+\frac{\gamma}{2}||\mathbf{x}\!-\!\mathbf{y}||^2,\label{eq:strongconvexity11}
\end{align}
\end{subequations}
for all $\mathbf{x},\mathbf{y}\in\mathbb{R}^n$ and  $\bth\in\Theta$.
In addition, there exists $M>0$ such that 
\begin{equation}\label{boundpartial}
 ||\frac{\partial}{\partial\bth} \nabla_{\bx} f(\mathbf{x},\bth) ||\leq M,
\end{equation}
for all $\mathbf{x}\in\mathbb{R}^n$ and all $\bth\in\Theta$.
\QEDB
\end{assumption}

\vspace{0.1cm}
The conditions \eqref{eq:globallipschitz1} and \eqref{eq:strongconvexity11} imply the smoothness and strong convexity of $f$ with respect to $\bx$, respectively, uniformly on $\theta$. Since $\Theta$ is a compact set, the uniformity assumption is not restrictive since one could obtain $\ell$ (resp. $\gamma$) by maximizing (resp. minimizing) $\theta$-dependent Lipschitz constants (resp. strong convexity constants) over $\Theta$.  These conditions are commonly assumed in time-varying optimization problems and enable exponential practical input-to-state stability bounds for the trajectories of the P-GZO dynamics \eqref{eq:pzo}. The following theorem states the tracking performance of the P-GZO dynamics \eqref{eq:pzo}, while preserving the Practical Safety property \eqref{eq:safetybound}. The proof is presented in Section \ref{sec:proof:track}.

%The smoothness condition implied by \eqref{globallipschitz1} and the strong convexity property \eqref{strongconvexity11} are common in time-varying model-based optimization problems, and they will allow us to establish a ``safe'' exponential practical input-to-state stability result for the model-free dynamics \eqref{algo1}, where the trajectories $t\mapsto\mathbf{x}(t)$ will now converge to a neighborhood of $\mathbf{x}^*(\theta)\in\mathcal{X}$. The size of the neighborhood will depend on the (point-wise) norm of $\dot{\mathbf{x}}^*(\theta)$, thus recoverging \eqref{KLbound1} when $\mathbf{x}^*$ is constant. This result is presented in the following theorem.
%

\vspace{0.1cm}
\begin{theorem} \label{thm:track} % \label{theorem2}
Consider the system dynamics \eqref{eq:pzo} and \eqref{eq:exosystem} with the flow set  
$\mathbf{C}_1\times\Theta$. Suppose that Assumptions \ref{ass:con_sm}-\ref{ass:track} hold.  Then, there exists $c>0$ such that for any $\Delta>\nu>0$, there exists $\hat{\varepsilon}_{\xi}>0$ such that for all $\varepsilon_{\xi}\in(0,\hat{\varepsilon}_{\xi})$, there exists $\hat{\varepsilon}_a>0$ such that for all $\varepsilon_a\in(0,\hat{\varepsilon}_a)$, there exists $\hat{\varepsilon}_{\omega}>0$ such that for all $\varepsilon_{\omega}\in(0,\hat{\varepsilon}_{\omega})$,
 every solution $\mathbf{z}$ of the P-GZO dynamics with $\mathbf{z}(0)\in \mathbf{C}_1\cap \big(( \mathbf{w}^*(0)+\Delta\mathbb{B})\times\mathbb{T}^n\big)$ is complete and satisfies the Practical Safety property \eqref{eq:safetybound}, and also: \vspace{-3pt}
\begin{align}
&\text{({Practical Tracking}):}\nonumber \\ &\qquad \limsup\limits_{t\rightarrow\infty} ||\mathbf{x}(t)-\mathbf{x}^*(\bth(t))||\leq c\cdot \sup_{t\geq0}||\dot{\bth}(t)||+\nu. \label{eq:ISSbound1}
%\\
%&\text{({Practical Safety}):}\ \ \mathbf{x}(t)\!\in\!\mathcal{X},\ \hat{\mathbf{x}}(t)\!\in\!\mathcal{X}\!+\!\varepsilon_a\mathbb{B},\ \forall t\geq 0, \label{eq:safety2}
\end{align}
where  $\mathbf{w}^*(0):= \big(\mathbf{x}^*(\bth(0)),  \nabla f(\mathbf{x}^*(\bth(0))) \big)$.
%\begin{enumerate}
%\item Practical Tracking:
%\begin{equation}\label{ISSbound1}
%\lim\sup_{t\to\infty}|\mathbf{x}(t)-\mathbf{x}^*(\theta(t))|\leq c_3 \sup_{t\geq0}|\dot{\theta}(t)|+\nu.
%\end{equation}
%where $\mathbf{z}^*=\{\mathbf{x}^*\}\times\{\nabla f(\mathbf{x}^*)\}$.
%\item Practical Safety:
%\begin{equation}\label{safety2}
%\mathbf{x}(t)\in\mathcal{X},~~\hat{\mathbf{x}}(t)\in\mathcal{X}+\varepsilon_a\mathbb{B},
%\end{equation}
%for all $t\geq0$.  \end{enumerate}
\hfill  \QEDB
\end{theorem}

The proof of Theorem \ref{thm:track} relies on  
input-to-state stability (ISS) tools for perturbed systems, which have been recently exploited to study other model-free optimization problems, e.g., \cite{labar2022iss,suttner2022robustness,scheinker2012extremum}. Due to the local boundedness of $\Pi(\cdot)$ and the compactness of $\Theta$, the function $t\mapsto \dot{\bth}(t)$ is uniformly bounded,  and thus the term $\sup_{t\geq0}||\dot{\bth}(t)||$ in \eqref{eq:ISSbound1} is well-defined and bounded by $\varepsilon_{\theta}c\rho_{\theta}$, where $|\Pi(\Theta)|\subset \rho_{\theta}\mathbb{B}$.

%However, in the context of model-free problems with hard constraints, they seem to be mostly unexplored.
%
%
%
%
\vspace{0.1cm}

\begin{figure}[t!]
    \centering
\includegraphics[width=0.5\textwidth]{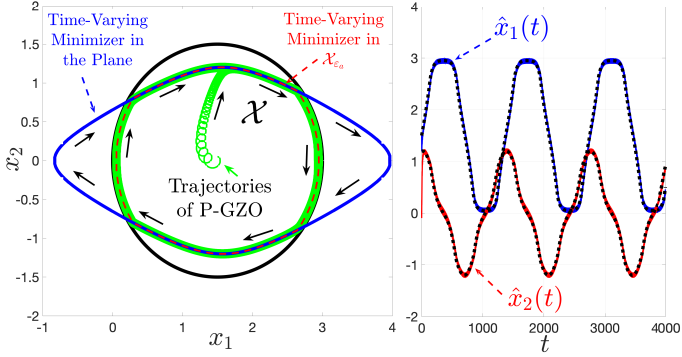}
      \caption{\small{Trajectories of P-GZO dynamics using a shrunk feasible set $\mathcal{X}_{\varepsilon_a}$, satisfying $\mathcal{X}_{\varepsilon_a}+\varepsilon_a\mathbb{B}\subset\mathcal{X}$. Left: The trajectories generated by the algorithm track the minimizer of $f$ inside the feasible set $\mathcal{X}$. Right: Evolution in time of the trajectories $x$. The optimal trajectories are shown with dotted lines.}}
    \label{fig:track}
    \vspace{-0.3cm}
\end{figure}

\begin{example} \label{example:pzo:track}
To illustrate the tracking performance of the P-GZO dynamics, we consider a simple problem in the plane, where $f(\mathbf{x},\bth)=({x}_1-{x}_{1}^*(\bth))^2+({x}_2-{x}_{2}^*(\bth))^2$ and the feasible set is the disk $\sX:=\{\bx\in\mathbb{R}^2: (x_1-1.5)^2+x_2^2\leq \frac{9}{4}\}$. Let $\bth$ be generated by the dynamics 
 $\dot{\theta}_1=\varepsilon_{\theta}\sin(2\theta_2)$, $\dot{\theta}_2=\frac{\varepsilon_{\theta}}{2}\cos(\theta_1)$, with $\varepsilon_{\theta}\!=\!1\times 10^{-2}$ and let  $x_i^*(\bth):=\theta_i$, $i\in[2]$.
 %In Figure \ref{fig:Fig1}, we used a high-frequency $\omega=2\pi\kappa_i/\varepsilon_{\omega}$ and a small amplitude $\varepsilon_a\in \mathcal{O}(10^{-3})$, so that the trajectories $\mathbf{x}$ approximate those of \eqref{projectedgradientflow1}, shown in dotted lines. In Figure \ref{fig:Fig2},
 To ensure strict safety, we use an $\varepsilon_a$-shrunk feasible set $\sX_{\varepsilon_a}:= \{\bx: (x_1-1.5)^2+x_2^2\leq (\frac{3}{2}-\varepsilon_a)^2\}$.
 Figure \ref{fig:track} shows the time-varying optimizer trajectory $t\mapsto\bx^*(t)$ and the solution trajectory $t\mapsto\bx(t)$ of the P-GZO dynamics \eqref{eq:pzo} under frequencies that are not necessarily too large, e.g., $\varepsilon_{\omega}\sim\mathcal{O}(10)$, and moderately small amplitudes, e.g., $\varepsilon_a\sim \mathcal{O}(10^{-2})$, which is a situation that is common in practical applications with computational limitations. The right plot shows the trajectories  $t\mapsto\hat{\bf{x}}(t)$, and it can be observed that it closely tracks $\bx^*(t)$ in the interior and the boundary of $\mathcal{X}$.  \QEDB
\end{example}

\vspace{0.1cm}
\begin{remark}
The bound \eqref{eq:ISSbound1} highlights the role of the rate of change of $\theta$ on the tracking error. When $\theta$ changes rapidly the P-GZO algorithm will generate a larger residual tracking error. On the other hand, as $\varepsilon_{\theta}\to0$ in \eqref{eq:exosystem}, such error will vanish, leading only to the residual bound $\nu$, which can be made arbitrarily small by decreasing $\varepsilon_{\omega},\varepsilon_a$. Note that decreasing $\varepsilon_{\theta}$ is equivalent to increasing the gains $k_x$ and $\frac{1}{\varepsilon_{\xi}}$ in the algorithm after a suitable change of time scale.  \QEDB 
\end{remark}
 
%We show the behavior of the algorithm under frequencies that are not necessarily too large, e.g., $\varepsilon_{\omega}\sim \mathcal{O}(10^0)$, and amplitudes that are not necessarily too small $\varepsilon_a\sim \mathcal{O}(10^{-1})$, which is a situation that is common in practical applications with computational limitations. To guarantee safety, the P-GZO algorithm implements a shrunk feasible set $\hat{\mathcal{X}}$ that satisfies $\hat{\mathcal{X}}+\varepsilon_a\mathbb{B}\subset\mathbb{R}^n$. As observed, the algorithm approximately tracks the desired optimal point $t\mapsto \mathbf{x}^*(\theta(t))$ in the interior and in the boundary of $\mathcal{X}$.
%
%\vspace{0.1cm}
\subsection{Switching Objective Functions}
This subsection considers the problem setting with switching objective functions. Depending on the information available to the decision-maker,  the objective function in \eqref{eq:gen:reduce} is drawn from a finite collection of functions $\{f_q(\mathbf{x})\}_{q\in Q}$. The selection of the current function to be optimized at each time $t$ might be performed by an external entity, leading to passive switching, or by the decision-maker, leading to active switching. In both cases, we show that provided the minimizers and critical points coincide across functions, the P-GZO dynamics can achieve safe optimization in a model-free way. % Real-time optimization problems with switching costs and safety constraints emerge in many engineering problems. For example, in the economic dispatch problem in electric power systems, the generation fuel costs and electricity prices can change over time. %, but sudden small price changes  may not lead to   different optimal dispatch solutions. 

\begin{figure}[t!]
    \centering
\includegraphics[width=0.43\textwidth]{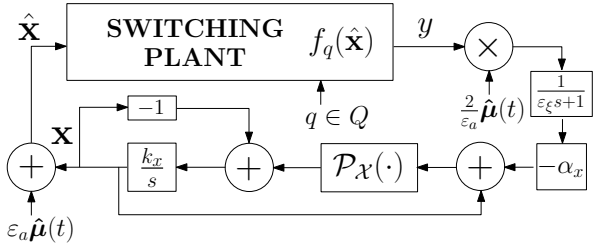}
      \caption{\small{Scheme of P-GZO dynamics with switching objectives.}}
    \label{fig:FigSwitch}
    \vspace{-0.3cm}
\end{figure}

Under switching objective functions,  the dynamics of the low-pass filter \eqref{eq:pzo:xi} become
%Under the smoothness and strong convexity properties \eqref{globallipschitz1}-\eqref{strongconvexity11}, the P-GZO dynamics \eqref{algo1} can also be implemented in applications where, depending on the available information to the user, a finite family of cost functions $\{f_q(\mathbf{x})\}_{q\in Q}$, sharing the same minimum point $\mathbf{x}^*$, are used by the ZO dynamics to find the point $\mathbf{x}^*$. In this case, the dynamics of the low-pass filter \eqref{filter1} become:
%
\begin{align}\label{eq:xi:track}
\dot{\bxi}=\frac{1}{\varepsilon_{\xi}}\Big(-\bxi+\frac{2}{\varepsilon_a}f_q(\hat{\mathbf{x}})\hat{\bmu}\Big),
\end{align}
where $q$ is a switching signal that selects from the set of indices $Q:=\{1,2,\ldots,\bar{q}\}$, with $\bar{q}<\infty$, the function $f_q$ to be used in the model-free algorithm at each time $t$, see Figure \ref{fig:FigSwitch}.
%where $q$ is a switching signal that selects from the index set  $Q:=\{1,2,\ldots,\bar{q}\}$ with $\bar{q}<\infty$. 
%the function $f_q$ to be used in the model-free algorithm at each time $t$, see Figure \ref{fig:FigSwitch}. 
This switching signal is generated by the following hybrid dynamical system \cite{Goebel:12}:
\begin{subequations}\label{eq:hybridautomaton}
\begin{align}
&(q,\tau)\in Q\times[0,N_0],~~~\dot{q}=0,~~~\qquad\ \dot{\tau}\in\left[0,\frac{1}{\tau_d}\right],\\
&(q,\tau)\in Q\times[1,N_0],~~~q^+\in Q\backslash
\{q\},~~\tau^+=\tau-1, \label{timerswitching}
\end{align}
\end{subequations}
where the state $\tau$ is a timer indicating when the signal $q$ is allowed to switch via \eqref{timerswitching}. In \eqref{eq:hybridautomaton}, $\tau_d>0$ is called the dwell-time, and $N_0\in\mathbb{Z}_{\geq1}$ is the chatter bound. As shown in \cite[Ch.2]{Goebel:12}, the hybrid system \eqref{eq:hybridautomaton} guarantees that every switching signal $q$ satisfies an average dwell-time (ADT) constraint. In particular, for every pair of times $(t_1,t_2)$ with $t_2>t_1$, every solution of \eqref{eq:hybridautomaton} satisfies:
\begin{equation}
\mathcal{S}(t_1,t_2)\leq \frac{1}{\tau_d}(t_2-t_1)+N_0, 
\end{equation}
where $\mathcal{S}(t_1,t_2)$ is the number of switches between times $t_1$ and $t_2$. The following theorem establishes the convergence and safety properties of the P-GZO dynamics under switching objectives. For simplicity, we consider the static optimization case when the optimizer $\bx^*$ is not time-varying but remains the same, and we omit the dependence of $\bx$ on discrete-time indices, which is typical in hybrid systems of the form \eqref{eq:hybridautomaton}. The proof of Theorem \ref{thm:switchedprojected} is provided in   Section \ref{sec:proof:switch}.
%Under these switching signals, we can establish the following result for the P-GZO algorithm whenever the critical point is also common. For simplicity, we present the case when $\Theta$ is a singleton and $\varepsilon_0=0$, i.e., the optimizer is not time-varying.
%
\vspace{0.1cm}
\begin{theorem}\label{thm:switchedprojected}
Consider the system dynamics \eqref{eq:pzo:x}, \eqref{eq:pzo:mu}, \eqref{eq:xi:track}, and \eqref{eq:hybridautomaton}.
Suppose that all functions in $\{f_q(\bx)\}_{q\in Q}$ are strongly convex and smooth,  Assumptions \ref{ass:con_sm} and \ref{ass:finite} hold for each of them, and they share
\begin{enumerate}
\item [(a)] common minimizer: ${\bx^*}=\arg\min_{x\in\mathcal{X}} f_q(\bx)$,~for all indices $q\in Q$;
\item [(b)] common critical point: $\bxi^*=\nabla f_q(\bx^*)$, ~for all indices $q\in Q$.
\end{enumerate}
Then, for any $\Delta>\nu>0$, there exists $\hat{\varepsilon}_{\xi}>0$ such that for all $\varepsilon_{\xi}\in(0,\hat{\varepsilon}_{\xi})$, there exists $\hat{\varepsilon}_a>0$ such that for all $\varepsilon_a\in(0,\hat{\varepsilon}_a)$, there exists $\hat{\varepsilon}_{\omega}>0$ such that for all $\varepsilon_{\omega}\in(0,\hat{\varepsilon}_{\omega})$,
 every solution $\bz(t)$ of the P-GZO dynamics \eqref{eq:pzo:x}, \eqref{eq:pzo:mu}, \eqref{eq:xi:track} with $\mathbf{z}(0)\in  \mathbf{C}_1 \cap \big(((\bx^*,\bxi^*)+\Delta\mathbb{B})\times\mathbb{T}^n\big)$ is complete and satisfies the Practical Safety property \eqref{eq:safetybound}, and also:
\vspace{-3pt}
\begin{align}
\text{({Practical Stability under ADT Switching}):}\nonumber \\  \limsup\limits_{t\rightarrow\infty} ||\mathbf{x}(t)-\mathbf{x}^*||\leq \nu.\qquad\quad \label{eq:switch} \hfill\QEDB
\end{align}  
%\begin{enumerate}
%\item Practical Convergence Under ADT Switching:
%\begin{equation*}
%\lim\sup_{t\to\infty}~|\mathbf{x}(t)-\mathbf{x}^*|\leq \nu,
%\end{equation*}
% where $\mathbf{z}^*=\{\mathbf{x}^*\}\times\{\xi^*\}$.
% \item Practical Safety:
%\begin{equation*}
%\mathbf{x}(t)\in\mathcal{X}~~~\text{and}~~~\hat{\mathbf{x}}(t)\in\mathcal{X}+\varepsilon_a\mathbb{B} 
%\end{equation*}
%
%for all $t\geq0$. 
% \end{enumerate} 
\end{theorem}

\vspace{0.1cm}
For \emph{unconstrained} optimization problems, it has been shown in \cite[Section 5.2]{PovedaTAC17B} and \cite{galarza2021extremum} that switched  ES algorithms are stable when each mode is stable and the switching is sufficiently slow. The novelty of Theorem \ref{thm:switchedprojected} lies in the incorporation of constraints  into the switching zeroth-order dynamics via projection maps.  Moreover, as shown in the analysis, the rate of convergence in \eqref{eq:ISSbound1} and \eqref{eq:switch} is actually exponential.
%
%\begin{remark}\label{remarkexample}
%

 Real-time optimization problems with switching costs and safety constraints emerge in many engineering problems. For example, in the economic dispatch problem in electric power systems, the generation fuel costs and electricity prices can change over time leading to changes in the landscape of the cost functions, but sudden small price changes may not lead to different optimal dispatch solutions. When the equilibrium points $x_q^*$ are distinct for each $q\in Q$, but they are all confined to a small $\delta$-neighborhood, assumptions (a) and (b) can be relaxed at the expense of obtaining a semi-global practical result with respect to $\delta$.
%
%\end{remark}
%
%Figure \ref{} shows a numerical example where the costs $f_q$ are given by .
%
\vspace{-0.3cm}
\subsection{Projected Primal-Dual ZO Dynamics with Lipschitz Projections} \label{sec:pdzo}
We now consider the complete optimization problem \eqref{eq:gen}, including the inequality constraints \eqref{eq:gen:ineq}. To solve  this problem, we first introduce the dual variable $\bla:=(\lambda_j)_{j\in[m]} \in\R_+^m$ for the inequality constraints \eqref{eq:gen:ineq},
and we formulate the  saddle 
 point problem  \eqref{eq:saddle}:
\begin{align} \label{eq:saddle}
     \max_{\bla\in \R_+^m}\,  \min_{\bx\in\mathcal{X}}\, L(\bx,\bla):=f(\bx)+  \bla^\top \bg(\bx),
\end{align}
where  $L(\bx,\bla)$ is the  Lagrangian function. Denote $\by:= [\bx;\bla]$, 
define $\mathcal{Y}:= \mathcal{X}\times \R_+^{m}$ as the feasible set  of $\by$, and denote $\sY^*$ as the set of the saddle points that solve  (\ref{eq:saddle}). 
By strong duality (implied by Assumptions \ref{ass:con_sm} and \ref{ass:finite}), the $\bx$-component of any saddle point $\by^*:= [\bx^*;\bla^*]\in\sY^*$ of \eqref{eq:saddle} is an optimal solution to problem \eqref{eq:gen}.

Similar to the study of the P-GZO dynamics \eqref{eq:pzo} in Section \ref{sec:pzodesign}, 
we now consider the \emph{projected primal-dual zeroth-order} (P-PDZO) dynamics \eqref{eq:c:pdzd} for the solution of problem \eqref{eq:gen}:
\begin{subequations} \label{eq:c:pdzd}
\begin{align}
        \dot{\bx} &= k_x\Big[\sP_{\sX\ }\big( \bx - \alpha_x \,\bxi_1 \big) -\bx\Big], \label{eq:c:pdzd:x} 
  \\
   \dot{\bla} & =k_\lambda\Big[ \sP_{\R_+^m\!}\big( \bla + \alpha_\lambda \bxi_2
     \big) -\bla \Big],\label{eq:c:pdzd:la}\\
          \dot{\bxi}_1& =\! \frac{1}{\varepsilon_\xi} \Big[ \!-\bxi_1 + \frac{2}{\varepsilon_a} \big(f(\hat{\bx}) + 
          \bla^\top \bg (\hat{\bx}) \big)\hat{\bmu}    \Big],
    \label{eq:c:pdzd:psi}\\ 
 \dot{\bxi}_2 &=\! \frac{1}{\varepsilon_\xi}\Big[\! -\bxi_2 + \bg (\hat{\bx}) \Big],  \label{eq:c:pdzd:psi2}
 \\
 \dot{ {\bmu}} &=\frac{1}{\varepsilon_{\omega}}\Lambda_{\kappa} {\bmu},\label{eq:c:pdzd:mu}
\end{align}
\end{subequations}
where  the parameters are defined in the same way as \eqref{eq:pzo}, and $\hat{\bx}$ and $\hat{\bmu}$ are defined as \eqref{eq:hatx} and \eqref{eq:hatmu}, respectively. Thus, the P-PDZO dynamics  \eqref{eq:c:pdzd} is restricted to evolve in the flow set
\begin{equation}\label{flowset2}
\mathbf{C}_2=\mathcal{X}\times\mathbb{R}_+^m\times\mathbb{R}^n\times\mathbb{R}^m\times\mathbb{T}^n. 
\end{equation}

The P-PDZO dynamics \eqref{eq:c:pdzd} can be regarded as a generalization of the P-GZO dynamics \eqref{eq:pzo} that further incorporates the inequality constraint \eqref{eq:gen:ineq}. Hence, the properties of P-GZO are generally applicable to P-PDZO, such as those mentioned in Remarks \ref{rem:lowpass} and \ref{rem:safety}.  
The following lemma states the forward invariance of $\mathbf{C}_2$, which
directly follows by Lemma \ref{lem:forward} by replacing $\mathcal{X}$ with
  $\mathcal{X}\times\mathbb{R}_+^m$.

\vspace{0.1cm}
\begin{lemma}\label{lemmaPD1}
Suppose that Assumption \ref{ass:con_sm} holds. Let $\mathbf{z}:=(\mathbf{x},\bla,{\bxi}_1,{\bxi}_2,{\bmu})$ be a solution of the P-PDZO dynamics \eqref{eq:c:pdzd}. Then, $\mathbf{z}(t)\in \mathbf{C}_2$ and $\hat{\mathbf{x}}(t)\in\mathcal{X}+\varepsilon_a\mathbb{B}$ for all $t\in\text{dom}(\mathbf{z})$.  \QEDB
\end{lemma}
\vspace{0.1cm}

%
% \begin{figure}[t!]
%     \centering
%       \includegraphics[scale=0.45]{feed2.png}
%     \caption{Feedback interconnection between plant $(f,g)$ and the proposed P-PDZD   method \eqref{eq:c:pdzd}.}
%     \label{fig:feedback}
% \end{figure}
%
%Let  $\hat{\sZ}^*$ be the saddle point set that solves \eqref{eq:saddle00}. When $f$ and $g$ are convex, by strong duality (implied by Assumptions \ref{ass:con_sm} and \ref{ass:finite}), for any point $\mathbf{z}^*=(\mathbf{x}^*,\lambda^*)\in\hat{\sZ}^*$, the component $\mathbf{x}^*$ is a solution to problem \eqref{eq:gen}. Therefore, we will study the stability and convergence properties of the state $\mathbf{x}$ in \eqref{eq:c:pdzd2} with respect to this set. 

We study the stability properties of the P-PDZO dynamics \eqref{eq:c:pdzd} based on the stability of the nominal target system:
\begin{subequations} \label{eq:cppdgd}
\begin{align} 
         \dot{\mathbf{p}}_1 &= k_1\Big[\sP_{\sX}\Big( \mathbf{p}_1 - \alpha_{1}\big( \nabla f(\mathbf{p}_1) +   \nabla \bg(\mathbf{p}_1)^\top \bp_2 \big)   \Big) -\mathbf{p}_1 \Big], \label{eq:cppdgd:x} 
  \\
   \dot{\mathbf{p}}_{2} &  = k_2 \Big[\sP_{\R_+}^m\Big( \mathbf{p}_{2} + \alpha_{2}\, \bg(\mathbf{p}_1)
     \Big) -\mathbf{p}_{2}\Big], \label{eq:cppdgd:la}
\end{align}
\end{subequations}
where $\nabla \bg:= [\nabla^\top g_1;\cdots;\nabla^\top g_m ]$ is the Jacobian matrix. 
%$\bp_1\in \R^n$ and $\bp_2\in\R^m$ in  \eqref{eq:cppdgd} correspond to the primal state $\bx$ and dual state $\bla$ in \eqref{eq:c:pdzd}, respectively.  
The nominal system \eqref{eq:cppdgd} is 
a well-known projected saddle flow that has been widely studied in the literature \cite{cherukuri2016asymptotic,goebel2017stability}.

The following theorem shows that the component $\by$ of the solution of \eqref{eq:c:pdzd} will converge to  a neighborhood of the saddle-point set $\sY^*$ using only zeroth-order information of $f$ and $\bg$, provided $\sY^*$ is compact and uniformly globally asymptotically stable (UGAS) under the nominal system \eqref{eq:cppdgd}. The proof is presented in Section \ref{app:thm:escp}.

\vspace{0.1cm}
\begin{theorem}\label{thm:p-pdzo}
Let $\bp:=[\bp_1;\bp_2]$, and suppose that Assumptions \ref{ass:con_sm}-\ref{ass:frequenciesdithers} hold, and:  
\begin{enumerate}[(a)]
\item The saddle point set $\sY^*$ is compact;
\item Every solution of \eqref{eq:cppdgd} with $\bp(0)\in \sY$ is complete; 
\item System \eqref{eq:cppdgd} renders the saddle point set $\sY^*$ forward invariant and uniformly attractive.
\end{enumerate}
Then, for any $\Delta>\nu>0$, there exists $\hat{\varepsilon}_{\xi}>0$ such that for all $\varepsilon_{\xi}\in(0,\hat{\varepsilon}_{\xi})$, there exists $\hat{\varepsilon}_a>0$ such that for all $\varepsilon_a\in(0,\hat{\varepsilon}_a)$, there exists $\hat{\varepsilon}_{\omega}>0$ such that for all $\varepsilon_{\omega}\in(0,\hat{\varepsilon}_{\omega})$,
 every solution $\mathbf{z}(t)$ of the P-PDZO dynamics \eqref{eq:c:pdzd} with $\mathbf{z}(0)\in \mathbf{C}_2\cap ( \left(\mathcal{W}^*_2+\Delta\mathbb{B}\right)\times\mathbb{T}^n)$ is complete and satisfies:
\begin{align}
&\text{({Practical Stability}):}\quad \limsup\limits_{t\rightarrow\infty}  ||\mathbf{y}(t)||_{\mathcal{Y}^*}\leq \nu,\qquad \label{eq:pd:KLbound1}\\
&\text{({Practical Safety}):}\ \ \mathbf{y}(t)\!\in\!\mathcal{Y},\ \hat{\mathbf{x}}(t)\!\in\!\mathcal{X}\!+\!\varepsilon_a\mathbb{B},\ \forall t\geq 0, \label{eq:pd:safetybound}
\end{align}
where $\mathcal{W}_2^*:=\{(\mathbf{y},\bxi_1,\bxi_2)\in\R^{2(n+m)}:\mathbf{y}\in\mathcal{Y}^*,\,\bxi_1=\nabla f(\mathbf{x}) +   \nabla \bg(\mathbf{x})^\top \bla,\,\bxi_2 = \bg(\bx)\}.$
 \QEDB

\end{theorem}
% 

%From the result of Theorem \ref{theorem3}, it is evident that the variable $\mathbf{x}_2$ plays the role of dual variable in the model-free dynamics, while $\mathbf{x}_1$ plays the role of a primal variable. However, it is important to note that Theorem \ref{theorem3} does not guarantee convergence to a neighborhood of a saddle \emph{point}, but rather convergence to a neighborhood of the saddle \emph{set}. To establish semi-global practical stability, in a point-wise sense, one would need to use more specialized tools \cite{goebel2017stability}. Such tools do not seem to be fully developed yet for multi-time scale dynamical systems.
%

\vspace{0.1cm}
\begin{remark}
The assumption of having a compact saddle point set \( \sY^* \) in Theorem \ref{thm:p-pdzo} is common when employing singular perturbation or averaging techniques. For many practical applications, the feasible set \( \sX \) represents physical capacity limits or control saturation bounds, and is therefore naturally compact. In some cases, we can 
substitute the feasible region $\R_+^m$ of the dual state $\bla$ by the feasible box set  $[0,M_\lambda]^m$ with a sufficiently large $M_\lambda$ to encompass any solution of practical interest. \QEDB 
\end{remark}
\vspace{2pt}

As discussed in Section II-B, the vanilla ES algorithm \eqref{eq:simes} emulates the behavior of an $\mathcal{O}(\varepsilon_a)$-perturbed gradient flow. Similarly, the P-GZO dynamics \eqref{eq:pzo} emulate the behavior of an $\mathcal{O}(\varepsilon_a)$-perturbed projected gradient flow, and the P-PDZO dynamics \eqref{eq:c:pdzd} emulate the behavior of an $\mathcal{O}(\varepsilon_a)$-perturbed projected saddle flow. While model-based algorithms of this form have been extensively studied in the literature, continuous-time zeroth-order implementations of these dynamics with stability and safety guarantees were mostly unexplored. Since in many cases the stability properties of \eqref{eq:cppdgd} (see items (a)-(c) of Theorem 5) are established via the Krasovskii-LaSalle invariance principle, the result of Theorem 4 allows us to establish stability properties for the model-free algorithm with similar generality as their model-based counterparts.

\section{MODEL-FREE FEEDBACK OPTIMIZATION WITH DISCONTINUOUS PROJECTIONS} 
\label{sec:discontinuous}
In the previous section, all the ZO algorithms utilized the Euclidean projection onto the feasible set \( \mathcal{X} \), resulting in ordinary differential equations (ODEs) with Lipschitz continuous vector fields on the right-hand side. This continuity property facilitates the well-posedness and stability analysis of the ZO dynamics, since the existence and uniqueness of solutions are guaranteed by standard results for ODEs \cite[Theorem 3.1]{khalil2002nonlinear}. In this section, we now turn our attention to the study of another class of projected ZO dynamics that enforce the hard constraints \eqref{eq:gen} via \emph{discontinuous projection maps}. This type of projection has been extensively studied in the context of (discontinuous) model-based projected dynamical systems \cite{nagurney2012projected,hauswirth2021projected}.
 % We leverage recent results of projected gradient flows to show that similar model-free dynamics can be designed as in Section \ref{sec:algorithm}. 
 To simplify our presentation, we focus on problem \eqref{eq:gen:reduce}, which does not include the inequality constraints \eqref{eq:gen:ineq}. %As shown in Section \ref{sec:pdzo}, such inequality constraints can be handled via primal-dual dynamics or additional penalty terms in the objective $f$; see  \cite{8571158} for details.  
 %
% As mentioned in the previous section, such constraints can be handled via primal-dual dynamics or additional penalty terms on the cost $f$. They can also be directly incorporated into the set $\mathcal{X}$ under Assumptions \ref{ass:con_sm}-\ref{ass:finite}.
%
\subsection{GZO Dynamics with Discontinuous Projection}
\label{secdisconpro}
\begin{figure}[t!]
    \centering
\hspace{-0.5cm}\includegraphics[width=0.44\textwidth]{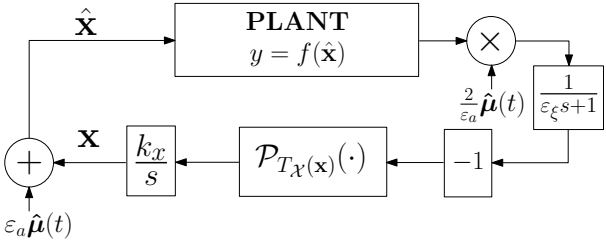}
      \caption{\small{Block diagram of the DP-GZO algorithm.}}
    \label{fig:eq:dpzo}
    \vspace{-0.3cm}
\end{figure}

To solve the reduced problem \eqref{eq:gen:reduce}, we consider the following ZO dynamics: 
\begin{subequations}\label{eq:dpzo}
\begin{align}
\dot{\mathbf{x}}&=k_x\mathcal{P}_{T_{\mathcal{X}}(\mathbf{x})}(-\,\bxi), \label{eq:dpzo:x} \\
\dot{ {\bxi}}&=\frac{1}{\varepsilon_{\xi}}\Big(-{\bxi}+\frac{2}{\varepsilon_a}f(\hat{\mathbf{x}}) \hat{\bmu}\Big),\label{eq:pzo:xi2} \\
\dot{ {\bmu}} &=\frac{1}{\varepsilon_{\omega}}\Lambda_{\kappa} {\bmu},\label{eq:pzo:mu2}
\end{align}
\end{subequations}
which are restricted to evolve in the flow set $\mathbf{C}_1$ defined in \eqref{flowset1}. In \eqref{eq:dpzo:x}, the mapping $\mathcal{P}_{T_{\mathcal{X}}(\mathbf{x})}(\cdot)$  projects the vector $-\bxi$ onto the tangent cone of the feasible set $\mathcal{X}$ at point $\bx$, i.e., $T_{\mathcal{X}}(\mathbf{x})$. As a result, the right-hand side of \eqref{eq:dpzo} is in general discontinuous, but it  guarantees that $\bx$ stays within the feasible set. Figure \ref{fig:eq:dpzo} shows a block diagram of the dynamics \eqref{eq:dpzo}. In conjunction with the flow set \eqref{flowset1}, we term the dynamics \eqref{eq:dpzo} as the \emph{discontinuous projected  gradient-based  zeroth-order} (DP-GZO) dynamics, and we study its stability and regularity properties using tools from differential inclusions and the following notion \cite[Definition 4.2]{Goebel:12}:
\vspace{0.1cm}
\begin{definition}
Consider the ODE $\dot{\mathbf{z}}=\bh(\mathbf{z})$, where $\mathbf{z}\in \mathbf{C}\subset \mathbb{R}^n$
 and $\bh:\mathbb{R}^n\to\mathbb{R}^n$ is locally bounded. The Krasovskii regularization of this ODE is the differential inclusion
\begin{equation}\label{krasovskiiregularization}
\mathbf{z}\in \overline{\mathbf{C}},~~\dot{\mathbf{z}}\in K(\mathbf{z}):= \bigcap_{\epsilon>0}\overline{\text{con}}~\bh(\left(\mathbf{z}+\epsilon\mathbb{B}\right)\cap \mathbf{C}),
\end{equation}
where, given a set $\sB$, $\text{con}(\sB)$ denotes its convex hull  and $\overline{\sB}$ denotes its closure. \QEDB
\end{definition}
\vspace{0.1cm}

The existence of solutions for the Krasovskii regularization of \eqref{eq:dpzo} is guaranteed by   well-posedness of \eqref{krasovskiiregularization} and standard viability results \cite[Theorems 3.3.4, 3.3.5]{viability1}. Moreover, it can be shown that system \eqref{krasovskiiregularization} accurately captures the limiting behavior of \eqref{eq:dpzo} under arbitrarily small additive perturbations on the states and dynamics via the so-called Hermes solutions \cite[Chapter 4]{Goebel:12}. This suggests that \eqref{krasovskiiregularization} provides a useful characterization of the solutions to the ODE $\dot{\bz}=\bh(\bz)$ under small perturbations, a setting that naturally emerges in the context of ZO dynamics. 

%\vspace{0.1cm}
Clearly, the solutions to \eqref{eq:dpzo} are also solutions of its Krasovskii regularization, but the converse is not always true. However, under mild regularity assumptions on \eqref{eq:dpzo}, it turns out that every solution of its Krasovskii regularization is also a standard (i.e., Caratheodory) solution of \eqref{eq:dpzo}, see \cite{hauswirth2021projected}. This fact allows us to study the behaviors of the DP-GZO dynamics \eqref{eq:dpzo} based on the following nominal ``target system'':
\begin{equation}\label{eq:projectedgradientflow}
\dot{\mathbf{p}}= k_x\cdot \mathcal{P}_{T_{\mathcal{X}}(\mathbf{p})}(-\nabla f(\mathbf{p})).
\end{equation}
The following theorem establishes the stability and (practical) safety properties for the DP-GZO dynamics \eqref{eq:dpzo}. The proof is presented in Section \ref{sec:proof:dpzo}.

\vspace{0.1cm}
\begin{theorem}  \label{thm:dpzo}  %\label{theorem4}
Suppose that Assumptions \ref{ass:con_sm}-\ref{ass:frequenciesdithers} hold and that $f$ is strictly convex. Then, for any $\Delta>\nu>0$, there exists $\hat{\varepsilon}_{\xi}>0$ such that for all $\varepsilon_{\xi}\in(0,\hat{\varepsilon}_{\xi})$, there exists $\hat{\varepsilon}_a>0$ such that for all $\varepsilon_a\in(0,\hat{\varepsilon}_a)$, there exists $\hat{\varepsilon}_{\omega}>0$ such that for all $\varepsilon_{\omega}\in(0,\hat{\varepsilon}_{\omega})$,
 every maximal solution $\mathbf{z}(t)$ of the DP-ZO dynamics \eqref{eq:dpzo}  with  $\mathbf{z}(0)\in \mathbf{C}_1\cap \left(\mathcal{W}_1^*+\Delta\mathbb{B}\right)\times\mathbb{T}^n)$  is complete and satisfies the Practical Convergence property \eqref{KLbound1} and the Practical Safety property \eqref{eq:safetybound}.  \hfill\QEDB
\end{theorem}
%

%\vspace{0.1cm}
%Although the results of Theorem \ref{thm:dpzo} for the DP-GZO dynamics \eqref{eq:dpzo} are similar to Theorem \ref{thm:pzo} for the P-GZO dynamics \eqref{eq:pzo}, their analyses and proof ideas are different, see  Section \ref{sec:proof:dpzo}.

\vspace{0.1cm}
\begin{remark}
% As $(\varepsilon_{\omega},\varepsilon_a,\varepsilon_{\xi})\to0^+$, the trajectories $\mathbf{x}$ generated by \eqref{disproj111} will approximate the trajectories of the projected gradient flow
% %
% \begin{equation}\label{projectedgradientflow}
% \dot{\mathbf{p}}=k_x \mathcal{P}_{T_{\mathcal{X}}(\mathbf{p})}(-\nabla f(\mathbf{p})),
% \end{equation}
% %
% which has been extensively studied in \cite{nagurney2012projected,hauswirth2021projected}. 
When $\varepsilon_{a},\varepsilon_{\xi},\varepsilon_{\omega}$ have small values, the $\mathbf{x}$-trajectories of the DP-GZO dynamics \eqref{eq:dpzo} emulate the trajectories of \eqref{eq:projectedgradientflow}. 
Since any closed and convex set $\mathcal{X}\subset\mathbb{R}^n$ is  Clarke regular \cite[Definition 2.2]{hauswirth2021projected} and prox-regular \cite[Definition 6.1]{hauswirth2021projected}, and since $f$ is locally Lipschitz under Assumption \ref{ass:con_sm}, the solutions to \eqref{eq:dpzo} and its Krasovskii regularization coincide and are unique \cite[Theorem 4.2]{thesisadrian}. Nevertheless, since uniqueness of solutions is not required in our analysis, the convexity of $\mathcal{X}$ and the strict convexity of $f$ could be relaxed to mere Clarke regularity and the assumption that all first-order critical points of \eqref{eq:gen} are optimal and also equilibria of \eqref{eq:projectedgradientflow}. \QEDB
\end{remark}
\begin{figure}[t!]
    \centering
\includegraphics[width=0.49\textwidth]{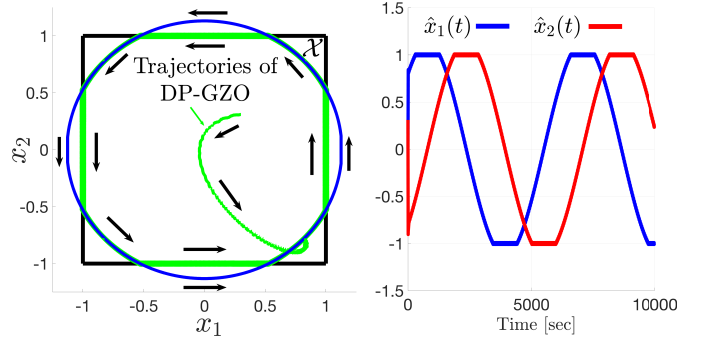}
      \caption{Illustration of the DP-GZO \eqref{eq:dpzo} restricted to a box $\mathcal{X}=[-1,1]\times[-1,1]$, with a slowly time-varying minimizer of $f$ corresponding to the blue circle trajectory.}
    \label{fig:eq:dpzo2}
    \vspace{-0.3cm}
\end{figure}
 \vspace{2pt}
\begin{example} \label{example:dpzo}
To illustrate the behavior of system \eqref{eq:dpzo}, we consider a simple example of problem \eqref{eq:gen:reduce}, where the feasible set $\sX:= [-1,1]\times [-1,1]$ is a box constraint and the objective function $f$ is the same as the one in Example \ref{example:pzo:track}, but now with $\dot{\theta}_1=-\varepsilon_p\theta_2$ and $\dot{\theta}_2=\varepsilon_p\theta_1$, $\varepsilon_p=1\times 10^{-3}$. As shown in  Figure \ref{fig:eq:dpzo2}, this exosystem makes the
 the minimizer of $f$ in $\mathbb{R}^n$  (i.e., $\arg\min_{\bx} f(\bx)$) a slowly varying signal that forms a circular trajectory, shown in blue color. It can be observed that the trajectory  $t\mapsto\mathbf{x}(t)$, shown in green, generated by the DP-GZO dynamics \eqref{eq:dpzo} tends to closely track the blue minimizer trajectory, but it stays in the feasible set at all times due to the projection map. Similar to Example \ref{example:pzo:track}, we have replaced $\sX$ with a shrunk set $\hat{\sX}_{\varepsilon_a}:= [-1+\varepsilon_a,1-\varepsilon_a]\times [-1+\varepsilon_a,1-\varepsilon_a]$ in  \eqref{eq:dpzo:x} to ensure the actual input $\hat{\bx}\in \sX$ all the time. Here, $\varepsilon_a=1\times 10^{-2}$, so the difference between $\hat{\sX}_{\varepsilon_a}$ and $\mathcal{X}$ is almost indistinguishable. The right plot   shows the trajectory of each of the components of $\hat{\bx}$.  \QEDB
\end{example}
 \vspace{2pt}
 
%Example \ref{example:dpzo} illustrates the convergence and  dynamic tracking performance of the DP-ZO dynamics method \eqref{eq:dpzo}. {\color{red}  The fact that system \eqref{eq:dpzo} is able to approximately track slowly varying bounded optimizers $\mathbf{x}^*$ follows directly by structural robustness properties of well-posed dynamical systems, discussed next.}

%in a problem \eqref{eq:gen} with box constraints and the same cost of Figure \ref{fig:Fig2}, where the optimal point that minimizes $f$ is slowly oscillating, generating a circular trajectory, shown in blue color. This optimal trajectory enters and leaves each of the sides of the feasible set. It can be observed that the trajectories approximately track the optimal solution of \eqref{eq:gen}. We also show in the right plot the trajectories of the actual input $\hat{\mathbf{x}}$, which, at the boundary of $\mathcal{X}$, performs small excursions with amplitude $\varepsilon_a$. Similar to the results of the previous section, for safety-critical applications the set $\mathcal{X}$ must be replaced by the largest set $\hat{\mathcal{X}}\subset \mathcal{X}$ such that $\hat{\mathcal{X}}+\varepsilon_a\mathbb{B}\subset\mathcal{X}$.

\vspace{-0.2cm}
\subsection{Structural Robustness}
The ZO  algorithms proposed in this paper rely heavily on function evaluations (or system output measurements) to steer the decision variable $\hat{\mathbf{x}}$ to an optimal solution of problem \eqref{eq:gen} or \eqref{eq:gen:reduce}.
Hence, suitable robustness properties are necessary  to handle small disturbances and noises that are inevitable in practice. The following result, i.e., Corollary \ref{coro:robust}, 
 indicates that, under the corresponding assumptions of Theorems \ref{thm:pzo}-\ref{thm:dpzo}, all the proposed ZO algorithms \eqref{eq:pzo}, \eqref{eq:c:pdzd}, \eqref{eq:dpzo}    are structurally   robust   to small bounded additive perturbations  on the states and dynamics. To state the corollary, we rewrite the ZO dynamics as a constrained ODE of the form $\mathbf{z} \in \mathbf{C},~ \dot{\mathbf{z}}=\bh(\mathbf{z} )$, and we consider their perturbed dynamics    \eqref{eq:perturbedODE}
 \begin{equation}\label{eq:perturbedODE}
\mathbf{z}+\mathbf{e}\in \mathbf{C},~~~\dot{\mathbf{z}}=\bh(\mathbf{z}+\mathbf{e})+\mathbf{e},
\end{equation}
 where $\bz$ is the state of the ZO dynamics, $\bh(\cdot)$ denotes the  vector field describing the right-hand side of the dynamics,  $\be$ is the additive noise, and $\mathbf{C}$ denotes the flow set.
\vspace{0.1cm}
\begin{corollary} \label{coro:robust}
    Under the assumptions and parameters of Theorems \ref{thm:pzo}-\ref{thm:dpzo}, there exists $\bar{e}>0$ such that for any  measurable function $\mathbf{e}(t):[0,+\infty)\!\to\! \R^{n}$ with $\sup_{t\geq 0} ||\mathbf{e}(t)||\leq \bar{e}$, the trajectory $\bz(t)$ of the perturbed ZO dynamics \eqref{eq:perturbedODE} satisfies the respective practical convergence bounds in Theorems \ref{thm:pzo}-\ref{thm:dpzo}. 
%\begin{align}\label{eq:robust}
%    ||\bz(t)-\hat{\bz}^*||\leq \beta(||\bz(0) -\hat{\bz}^*||,\,  t) +2\nu, \ \ \forall t\geq 0.
%\end{align}
% \begin{align}\label{eq:robust}
%     ||\bx(t)||_{\mathcal{A}}\leq \beta(||\bx(t_0)||_{\mathcal{A}},\,  t-t_0) +2\nu, \ \ \forall t\geq t_0.
% \end{align}
% %
% for all $t\geq t_0$. 
\QEDB
\end{corollary}

The result is a corollary of Theorems \ref{thm:pzo}-\ref{thm:dpzo} because the convergence,  well-posedness, and stability properties of the dynamics imply that, for each sufficiently large compact set of initial conditions $K$, and fixed parameters of the controller that induce the convergence bounds, there exists a compact set that is locally asymptotically stable under the nominal dynamics (the so-called Omega-limit set of $K$\footnote{See \cite[Definition 6.23]{Goebel:12} for the  notion of ``Omega-limit set of a set".}), and also semi-globally practical asymptotically stable as $\bar{e}\to0^+$ for the perturbed system \eqref{eq:perturbedODE} \cite[Chapter 7]{Goebel:12}. Similar robustness results have been studied in the literature of discontinuous systems \cite{shi2022finite}. However, we note that Corollary \ref{coro:robust} only shows the \emph{existence} of a sufficiently small $\bar{e}$, such that any additive disturbance bounded by $\bar{e}$ does not change drastically the convergence properties of the ZO algorithms. However, in practice, the explicit computation of this robustness bound is challenging and application-dependent. 

\vspace{-0.2cm}
\section{ANALYSIS AND PROOFS}
\label{section_analysis}
In this section, we present the proofs of our main results. Since the result of Theorem \ref{thm:pzo} can be seen as a particular case of Theorem \ref{thm:p-pdzo} when the set of inequality constraints \eqref{eq:gen:ineq} is empty, we first  present the proof of Theorem \ref{thm:p-pdzo}. Subsequently, we show how to adapt this proof to Theorem \ref{thm:pzo}. The proofs of Theorems \ref{thm:track}, \ref{thm:switchedprojected}, \ref{thm:dpzo} are based on the construction of suitable Lyapunov functions, and therefore are presented afterwards.

\vspace{-0.3cm}
\subsection{Proof of Theorem \ref{thm:p-pdzo}}
\label{app:thm:escp}
Let $\by:=[\bx;\bla]$, $\bxi:=[\bxi_1;\bxi_2]$, and $\bs:=[\by;\bxi]$. 
The P-PDZO dynamics (\ref{eq:c:pdzd:x})-\eqref{eq:c:pdzd:psi2} can be written in compact form as
\begin{align} \label{eq:original}
 \dot{\bs}=    \begin{bmatrix}
      \dot{\by} \\   \dot{\bxi}
    \end{bmatrix} = \begin{bmatrix}
    \mathbf{q}_1(\by,\bxi)\\
    \frac{1}{\varepsilon_{\xi}}(-\bxi + \mathbf{q}_2(\by,\bmu))
    \end{bmatrix}:= \mathbf{q}(\bs,\bmu),
\end{align}
where $\mathbf{q}_1$ captures the dynamics \eqref{eq:c:pdzd:x}-\eqref{eq:c:pdzd:la}, and $\bq_2$ is 
\begin{align}\label{eq:g2}
    \bq_2(\mathbf{y},\bmu):=\begin{bmatrix}
     \frac{2}{\varepsilon_a} \big(f(\hat{\bx})+\bla^\top \bg(\hat{\bx}) \big)\hat{\bmu} \\
  \bg(\hat{\bx}) 
    \end{bmatrix},
\end{align}
where $\bmu$ is generated by the oscillator \eqref{eq:c:pdzd:mu}, and $\hat{\mathbf{x}},\hat{\bmu}$ are defined in \eqref{eq:hatx} and \eqref{eq:hatmu}, respectively. We analyze the stability properties of this system using averaging and singular perturbation theory. We divide the analysis into the following three main steps.

\vspace{0.1cm}
\noindent \textbf{Step 1)} Let $\Delta>\nu>0$ and $\sY:=\mathcal{X}\times\mathbb{R}^m_+$, where without loss of generality we take $\nu<1$. Consider the compact set $[(\sY^*+\Delta \mathbb{B})\cap \sY]\times \Delta \B$ for the initial condition $\bs(0)$ %, where $\nu$ is the precision in Theorem \ref{thm:spas}
. Here, 
 $\sY^*+\Delta \mathbb{B}$ denotes 
 the union of all sets obtained by taking a closed ball of radius $\Delta$ around each point in $\sY^*$. %Basically, one can select sufficiently large $\Delta$ to cover all  feasible initial conditions. 

By items (a)-(c) in Theorem \ref{thm:p-pdzo}, the compact set $\mathcal{Y}^*$ is uniformly globally asymptotically stable
(UGAS) for the target system \eqref{eq:cppdgd} restricted to evolve in $\sY$, which is also a forward invariant set due to the projection mappings. Thus, there exists a class-$\mathcal{KL}$ function $\beta$ such that for any initial condition $\mathbf{p}(0)\in \sY$, the solutions $\mathbf{p}$ of   \eqref{eq:cppdgd} satisfy $||\mathbf{p}(t)||_{\sY^*}\leq \beta(||\mathbf{p}(0)||_{\sY^*},\,  t)$ for all $t\geq0$. Without loss of generality, let $\nu\in(0,1)$ and consider the set
\begin{align}\label{eq:Fset}
    \sF\!:=\! \Big\{\by\!\in \!\sY:   ||\mathbf{y}||_{\sY^*} \leq \beta\big( \max_{\bv\in \sY^*+\Delta \mathbb{B} } ||\bv||_{\sY^*},\, 0 \big)+1  \Big\}.
\end{align}
Note that the set $\sF$ is compact under the assumption that $\sY^*$ is compact.
Due to the boundedness of $\sF$, there exists $M_1>0$ such that $\sF\subset M_1 \B$. Let 
\begin{equation}\label{eqlll}
\bm{\ell}(\by)\!:=\!\begin{bmatrix}\nabla f(\bx) + \sum_{j=1}^m \lambda_j \nabla g_j(\bx) \\
   \bg(\bx)
    \end{bmatrix},
\end{equation}
and note that, by continuity of $\bm{\ell}$, there exists  $M_2>\max\{\Delta,1\}$ such that $||\bm{\ell}(\by)|| +1\leq M_2$ for all $||\by||\leq M_1$. Denote $M_3:= M_2+1$.
 We then study the behavior of system \eqref{eq:original} \emph{restricted to evolve in the compact set $\sF\times M_3\B$}.

\vspace{0.1cm}
\noindent \textbf{Step 2)} Since the solutions of the oscillator \eqref{eq:c:pdzd:mu} are given by \eqref{eq:periodicdither}, and $\mu_i(0)^2+\mu_{i+1}(0)^2=1$ for all $i\in\{1,3,\ldots,2n-1\}$,  system \eqref{eq:original} with small values of $\varepsilon_\omega$  is in standard form  for the application of averaging theory along the trajectories of $\bmu$.  The following Lemma \ref{lemma:aveg2} characterizes the average map of $\mathbf{q}_2$. The proof is presented in Appendix \ref{sec:app:lem}.
\begin{lemma}\label{lemma:aveg2}
The average of $t\mapsto\bq_2(\by,\bmu(t))$ is given by
\begin{align}
    \bar{\bq}_2(\by):=& \frac{1}{T}\!\int_0^T \!\bq_2(\mathbf{y},\bmu(t))\, dt
    = \bm{\ell}(\by) +\sO(\varepsilon_a),
\end{align}
where $\bm{\ell}$ is given by \eqref{eqlll}, and  $T>0$ is the common period of the dithers $\bmu$. \QEDB
\end{lemma}

\vspace{0.1cm}
Using Lemma \ref{lemma:aveg2},
we obtain the \emph{average dynamics} of \eqref{eq:original}:
\begin{align} \label{eq:realave}
    \dot{\bar{\bs}} = \!\begin{bmatrix} 
    \dot{\bar{\by}} \\
    \dot{\bar{\bxi}}
    \end{bmatrix} \!=\!\begin{bmatrix}
    \mathbf{q}_1(\bar{\by},\bar{\bxi})\\
    \frac{1}{\varepsilon_{\xi}}(-\bar{\bxi} + \bm{\ell}(\bar{\by}) +\sO(\varepsilon_a)  )
    \end{bmatrix},
\end{align}
where $\bar{\bs}:=[\bar{\by};\bar{\bxi}]$. We study \eqref{eq:realave}  restricted to evolve in the compact set $\sF\times M_3\B$. 
We treat the right-hand side of \eqref{eq:realave} as an $\sO(\varepsilon_a)$-perturbation of a nominal system with $\sO(\varepsilon_a)=0$. This nominal system is in the standard form for the application of singular perturbation theory \cite{teel2003unified}, with $\bar{\by}$ being the slow state, and $\bar{\bxi}$ being the fast state. The boundary layer dynamics of this nominal system, in the time scale $\tau=t/\varepsilon_{\xi}$, are 
\begin{align} \label{eq:linear}
    \frac{d \bar{\bxi}}{d \tau} = -\,\bar{\bxi} +\bm{\ell}(\bar{\by}),
\end{align} 
where $\bar{\by}$ is kept constant. This linear system \eqref{eq:linear} has a unique exponentially stable equilibrium point $\bar{\bxi}^*=\bm{\ell}(\bar{\by})$. As a result, the associated \emph{reduced system}  is derived as 
\begin{align}\label{eq:reduced}
    \dot{\bar{\by}} = \bq_1(\bar{\by},\, \bm{\ell}(\bar{\by})),
\end{align}
which is exactly the nominal target system \eqref{eq:cppdgd}. Under the assumptions of Theorem \ref{thm:p-pdzo}, system \eqref{eq:reduced} renders the set $\sY^*$ UGAS with $\beta\in\mathcal{K}\mathcal{L}$. By invoking stability results for singularly perturbed systems \cite[Theorem 2]{Wang:12_Automatica}, we can conclude that, as $\varepsilon_{\xi}\to 0^+$, the set $\sY^*\times {M}_3\mathbb{B}$ is semi-globally practically asymptotically stable (SGPAS) for the unperturbed average system \eqref{eq:realave} with $\sO(\varepsilon_a)=0$. Since system \eqref{eq:realave} has a continuous right-hand side , the perturbed average system  \eqref{eq:realave} also renders the set $\sY^*\times {M}_3\mathbb{B}$ SGPAS as $(\varepsilon_{\xi},\varepsilon_a)\to 0^+$, which is stated as Lemma \ref{lemma:y:sgpsa}.

\vspace{0.1cm}
\begin{lemma}\label{lemma:y:sgpsa} 
There exists $\beta\in\mathcal{K}\mathcal{L}$ such that for each $\nu>0$,
there exists $\hat{\varepsilon}_{\xi}>0$ such that for any $\varepsilon_{\xi}\in(0,\hat{\varepsilon}_{\xi})$, there exists $\hat{\varepsilon}_a>0$ such that for any $ \varepsilon_a\in (0,\hat{\varepsilon}_a)$, 
every solution $\bar{\bs}$ of the average system \eqref{eq:realave} (restricted in $\sF\times M_3\B$) with initial condition $\bar{\bs}(0)\in [(\sY^*\!+\!\Delta \mathbb{B})\cap \sY]\!\times\! \Delta \B$ satisfies 
\begin{align} \label{eq:y1con}
     ||\bar{\by}(t)||_{\sY^*}\!\leq \beta(||\bar{\by}(0)||_{\sY^*},\,  t) +\frac{\nu}{4},
\end{align}
for all $t\in \mathrm{dom}(\bar{\bs})$. \QEDB
\end{lemma}

\vspace{0.1cm}
Since the average system \eqref{eq:realave} is restricted in $\sF\times M_3\B$, we have $||\bar{\bxi}(t)||_{M_3\B} = 0$ for all $t\in \mathrm{dom}(\bar{\bs})$, which implies that $||\bar{\bs}(t)||_{\sY^*\times M_3\B} =  ||\bar{\by}(t)||_{\sY^*}$ for all $t\in \mathrm{dom}(\bar{\bs})$. Hence, \eqref{eq:y1con} implies that for all $ t\in \mathrm{dom}(\bar{\bs})$:
\begin{align*}
   ||\bar{\bs}(t)||_{\sY^*\times M_3\B}\leq \beta(||\bar{\bs}(0)||_{\sY^*\times M_3\B},\, t) +\frac{\nu}{4}.
\end{align*}
Next, we show the completeness of solutions of the average  system \eqref{eq:realave} by leveraging Lemma \ref{lemma:forward}, which follows as a special case of \cite[Lemma 5]{NunoShamma_2020}.
\begin{lemma} \label{lemma:forward}
Let $k,M_2>0$ be given and $\mathbf{u}:\R_+ \to M_2\B$ be a continuous function of time. Then, the set $M_2\B$ is forward invariant under the dynamics $\dot{\bxi} = k\cdot(-\bxi +\mathbf{u}(t))$. \QEDB
\end{lemma}

\vspace{0.1cm}
Under the initial condition $\bar{\bs}(0)\in [(\sY^*\!+\!\Delta \mathbb{B})\cap \sY]\!\times\! \Delta \B$, by \eqref{eq:y1con}, the trajectory $\bar{\mathbf{y}}$ of \eqref{eq:realave} satisfies $\bar{\mathbf{y}}(t)\in \mathrm{int}(\sF)$ for all $t\in \mathrm{dom}(\bar{\bs})$. This implies that  $||\bar{\mathbf{y}}(t)||\leq M_1$ and thus $||\bm{\ell}(\bar{\mathbf{y}}(t)) \!+\!\sO(\epsilon_a))||< M_2$ for all $t\in \mathrm{dom}(\bar{\bs})$, where, without loss of generality, we take $||\sO(\epsilon_a)||<1$ for all $\varepsilon_a\in (0,\hat{\varepsilon}_a)$. Using Lemma \ref{lemma:forward}, $\bar{\bxi}(t) \in M_2\B\subset \mathrm{int}(M_3\B)$ for all $t\geq 0$. Thus, under the given initialization, $\bar{\bs}$ satisfies 
\begin{align}\label{eq:unbs}
    \bar{\bs}(t)\in \mathrm{int}(\sF\times M_3\B), \quad \forall t\geq 0,
\end{align}
and thus it has an unbounded time domain. 

\vspace{0.2cm}
\noindent \textbf{Step 3)} Since  the set $\sY^*\times M_3\B$ is SGPAS for the average system \eqref{eq:realave} (restricted in $\sF\times M_3\B$) as $(\varepsilon_{\xi},\varepsilon_a)\to 0^+$, by using averaging theory for perturbed systems \cite[Theorem 7]{PovedaNaLi2019} it follows that for each pair of $(\varepsilon_{\xi},\varepsilon_a)$ inducing the bound \eqref{eq:y1con}, there exists $\hat{\varepsilon}_\omega>0$ such that for any $\varepsilon_\omega\in(0,\hat{\varepsilon}_\omega)$, the solution $\bs$ of the system \eqref{eq:original} (restricted to $\sF\times M_3\B$) satisfies 
\begin{align}
    ||\bs(t)||_{\sY^*\times M_3\B}   \leq \beta(||\bs(0)||_{\sY^*\times M_3\B},\, t) +\frac{\nu}{2},
\end{align}
for all $t\in \mathrm{dom}(\bs)$. Since  $ ||\bx(t)||_{\sY^*}= ||\bs(t)||_{\sY^*\times M_3\B}$ for all $t\in \mathrm{dom}(\bs)$, we obtain the bound  \eqref{eq:pd:KLbound1}.  The only task left is to show the completeness of solutions of the original  system \eqref{eq:original}. This can be done by using the following lemma, proved in Appendix \ref{app:omega}, as well as Lemma \ref{lemma:close}, which follows by \cite[Theorem 1]{Wang:12_Automatica}.
\vspace{0.1cm}
\begin{lemma}\label{lemma:omega}
There exists $\hat{\varepsilon}_{\xi}>0$ such that for any $\varepsilon_{\xi}\in(0,\hat{\varepsilon}_{\xi})$, there exists $\hat{\varepsilon}_a>0$ such that for any $ \varepsilon_a\in (0,\hat{\varepsilon}_a)$,  there exists a compact set $\Omega(\sF\times M_3\B)$ that is uniformly globally (pre)-asymptotically stable for the average system \eqref{eq:realave} restricted to $\sF\times M_3\B$, and which satisfies $\Omega(\sF\times M_3\B)\subset (\mathcal{Y}^*\times M_2\B) + \frac{\nu}{2} \B \subset \mathrm{int}(\sF\times M_3\B)$. \QEDB
\end{lemma}
\vspace{0.1cm}
\begin{lemma}\label{lemma:close}
Let $(\varepsilon_\xi,\varepsilon_a)>0$ take sufficiently small values such that Lemmas \ref{lemma:y:sgpsa} and \ref{lemma:omega} hold. Then, for each $\tau,\varepsilon>0$, there exists $\hat{\varepsilon}_\omega>0$ such that for all $\varepsilon_\omega\in (0,\hat{\varepsilon}_\omega)$  and for each solution $\bs$ to the original system \eqref{eq:original}, with $\mathbf{s}(0)\in [(\sY^*+\Delta\mathbb{B})\cap\sY]\times\Delta\mathbb{B}$, there exists a solution $\bar{\mathbf{s}}$ of the average system \eqref{eq:realave}, with $\bar{\mathbf{s}}(0)\in [(\sY^*+\Delta\mathbb{B})\cap\sY]\times\Delta\mathbb{B}$ such that $\mathbf{s}$ and $\bar{\mathbf{s}}$ are $(\tau,\varepsilon)$ close. %(restricted to $\sF\times M_3\B$) such that 
%
%\begin{align}\label{closeness1}
%    \sup_{t\in[0,\tilde{T}]}||\bx(t)-\bar{\bx}(t)||\leq \delta,~\sup_{t\in[0,\tilde{T}]}||\bxi(t)-\bar{\bxi}(t)||\leq \delta.
%\end{align}
%
%for all $t\in \mathrm{dom}(\bs)\cap \mathrm{dom}(\bar{\bs})$.
\hfill \QEDB
\end{lemma}
%The completeness of solution $\bs$ for the original system \eqref{eq:original} is guaranteed by taking $M_2$ sufficiently large.
%Thus Theorem \ref{thm:spas} is proved.
\vspace{0.1cm}

By Lemma \ref{lemma:omega}, there exists a $\mathcal{KL}$-class function $\tilde{\beta}$ such that every solution of the restricted average system \eqref{eq:realave} satisfies  $||\bar{\bs}(t)||_{\Omega( \sF\times M_3\B )}\leq \tilde{\beta}( ||\bar{\bs}(0)||_{\Omega( \sF\times M_3\B )},t)$, for all $t\in\mathrm{dom}(\bar{\bs})$. Therefore, by averaging theory \cite[Theorem 2]{Wang:12_Automatica}, there exists $\hat{\varepsilon}_\omega$ such that for all $\varepsilon_\omega\in(0,\hat{\varepsilon}_\omega)$, every solution of the original system \eqref{eq:original} with $\bs(0)\in [(\sY^*\!+\!\Delta \mathbb{B})\cap \sY]\!\times\! \Delta \B$, satisfies 
 \begin{align}\label{KLomegaaverage2}
    ||{\bs}(t)||_{\Omega( \sF\times M_3\B )}\leq \tilde{\beta}(||{\bs}(0)||_{\Omega( \sF\times M_3\B )},t)+\frac{\nu}{4},  
\end{align}
for all $t\in\text{dom}(\mathbf{s})$. Since by \eqref{eq:unbs} the trajectories $\bar{\bs}$ are complete if $\bar{\bs}(0)\in [(\sY^*\!+\!\Delta \mathbb{B})\cap \sY]\!\times\! \Delta \B$, using the closeness of solutions property of Lemma \ref{lemma:close} and the bound \eqref{KLomegaaverage2} it can be shown that, under the given initialization, $\bs$ satisfies $\bs(t) \in  \mathrm{int}(\sF\times M_3\B)$ for all $t\geq0$, i.e., every solution of the original system \eqref{eq:original} has an unbounded time domain. This establishes Theorem \ref{thm:p-pdzo}. 
% Thus, there exists time $T_3>0$ such that 
% %
% \begin{align}\label{eq:someb}
%     |{\bs}(t)|_{\Omega( \sF\times M_3\B )}\leq \frac{\nu}{2}, \ \forall  t\in [T_3,+\infty)\cap \mathrm{dom}({\bs}).
% \end{align}
% %
% Lemma \ref{lemma:close} indicates that all solutions of the original system \eqref{eq:original} remain $\delta$-close to some solution of the average system \eqref{eq:realave} on a compact time domain. Moreover, \eqref{eq:unbs} indicates that every solution of \eqref{eq:realave}, initialized in $[(\sY^*\!+\!\Delta \mathbb{B})\cap \sY]\!\times\! \Delta \B$, stays within $ \mathrm{int}(\sF\times M_3\B)$ for all $t\geq 0$. Thus, by applying Lemma \ref{lemma:close} with $\tilde{T} = T_3+1$, there exists $\hat{\varepsilon}_\omega>0$ such that for all $\varepsilon_\omega\in (0,\hat{\varepsilon}_\omega)$ we have
% %
% \begin{align}
%     \bs(t) \in  \mathrm{int}(\sF\times M_3\B), \quad \forall t\in [0,T_3+1].
% \end{align}
% %
% In addition, by \eqref{eq:someb}, we also have $ \bs(t) \in  \mathrm{int}(\sF\times M_3\B)$ for all $t\geq T_3$. Therefore, every solution $\bs(t)$ of the original system, under the given initialization, has an unbounded time domain. Thus Theorem \ref{theorem3} is proved.
%\hfill $\blacksquare$
%

\vspace{-0.3cm}
\subsection{Proof of Theorem \ref{thm:pzo}} \label{sec:pzo:proof}
The proof of Theorem \ref{thm:pzo} follows the same steps as the proof of Theorem \ref{thm:p-pdzo}. In this case, $\mathbf{q}_2:=\frac{2}{\varepsilon_a}f(\hat{\mathbf{x}})\hat{\bmu}$, and the function $\bm{\ell}(\cdot)$ in \eqref{eq:realave} becomes $\bm{\ell}(\mathbf{x})=\nabla f(\mathbf{x})$. Therefore, system \eqref{eq:reduced} is precisely the Lipschitz continuous projected gradient flow \eqref{eq:nominal:pzo}. Under conditions (a)-(b) in Theorem \ref{thm:pzo}, there exists a class $\mathcal{K}\mathcal{L}$ function $\beta$ such that Lemma \ref{lemma:y:sgpsa} holds. The rest of the proof follows the same steps as the proof of Theorem \ref{thm:p-pdzo}.

\vspace{-0.2cm}
\subsection{Proof of Theorem \ref{thm:track}} \label{sec:proof:track}
Following similar computations as in Lemma \ref{lemma:aveg2}, we compute the average dynamics of \eqref{eq:pzo}. In this case, we obtain
\begin{subequations}\label{auxsysthm22}
\begin{align}
\dot{\bar{\mathbf{x}}}&=k_x \mathcal{P}_{\mathcal{X}}\left(\mathbf{\bar{x}}-\alpha_x \bar{\bxi}\right)-k_x\mathbf{\bar{x}},\\
\dot{\mathbf{\bar{\bxi}}}&=\frac{1}{\varepsilon_{\xi}}\left(-{\bar{\bxi}}+\nabla_{\bar{\bx}} f(\mathbf{\bar{x}},\bar{\bth})+\mathcal{O}(\varepsilon_a)\right), \\
\dot{\bar{\bth}}&=\varepsilon_{\theta}\Pi(\bar{\bth}),
\end{align}
\end{subequations}
which evolve in the flow set $\mathbf{C}_3:=\mathcal{X}\times\mathbb{R}^n\times\Theta$, and which can be seen as an $\mathcal{O}(\varepsilon_a)$-perturbed two-time scale system with respect to the parameter $\varepsilon_{\xi}$. Next, we establish a key lemma for the average system \eqref{auxsysthm22}, which relies on singular perturbation tools \cite[Ch. 11]{khalil2002nonlinear}.

\vspace{0.1cm}
\begin{lemma}\label{auxlemma9}
For system \eqref{auxsysthm22} with $\mathcal{O}(\varepsilon_a)=0$ and $\varepsilon_{\theta}>0$, there exist $\hat{\alpha}_x, \hat{k}_x,\hat{\varepsilon}_{\xi}>0$, such that for all $\alpha_x\in(0,\hat{\alpha}_x)$, all $k_x\in(0,\hat{k}_x)$ and all $\varepsilon_{\xi}\in(0,\hat{\varepsilon}_{\xi})$, every solution satisfies
\begin{align*}
||{\bar{\bs}}(t)-{\bs}^*(\bar{\bth}(t))||\leq c_1||{\bar{\bs}}(0)&-{\bs}^*(\bar{\bth}(0))||e^{-c_2 t}\\
&\qquad +c\sup_{0\leq \tau\leq t}||\dot{\bar{\bth}}(\tau)||,
\end{align*}
where $\mathbf{\bar{s}}:=(\mathbf{\bar{x}},{\bar{\bxi}})$,  $\mathbf{s}^*(\bar{\bth}):=(\mathbf{x}^*(\bar{\bth}),\nabla_{\bx} f(\mathbf{x}^*(\bar{\bth}),\bar{\bth}))$, and $c_1,c_2,c>0$. \QEDB
\end{lemma}

\vspace{0.1cm}\noindent 
\textbf{Proof:} We introduce the error variables $\tilde{\mathbf{x}}:=\mathbf{\bar{x}}-\mathbf{x}^*(\bar{\bth})$ and $\tilde{\bxi}:=\mathbf{\bar{\bxi}}-\nabla_{\bx}f(\tilde{\mathbf{x}}+\mathbf{x}^*,\mathbf{\bar{\bth}})$, which leads to the error dynamics
\begin{subequations}\label{errordynamicsaux}
\begin{align}
\dot{\tilde{\mathbf{x}}}&=k_x\mathcal{P}_{\mathcal{X}}\left(\tilde{\mathbf{x}}+\mathbf{x}^*-\alpha_x(\tilde{\bxi}+\nabla_{\bx} f(\mathbf{\tilde{x}}+\mathbf{x}^*,\bar{\bth}))\right)\notag\\
&~~~-k_x (\tilde{\mathbf{x}}+\mathbf{x}^*)-\dot{\mathbf{x}}^*, \label{eq:error:x} \\
\dot{\tilde{\bxi}}&=-\frac{1}{\varepsilon_{\xi}}\tilde{\bxi}-\frac{d}{dt}\nabla_{\bx} f(\tilde{\mathbf{x}}+\mathbf{x}^*,\bar{\bth}).
\end{align} \label{eq:error:xi}
\end{subequations}
We study the stability properties of \eqref{errordynamicsaux} with respect to the origin. To this end, we consider the composite Lyapunov function: 
\begin{align} \label{eq:W:lya}
    W(\tilde{\mathbf{x}},\tilde{\bxi})=(1-\lambda)V_x(\tilde{\mathbf{x}})+\lambda V_{\xi}(\tilde{\bxi}),
\end{align}
where $\lambda\in(0,1)$,  $V_x(\tilde{\mathbf{x}})=\frac{1}{2}||\tilde{\mathbf{x}}|||^2$ and $V_{\xi}(\tilde{\bxi})=\frac{1}{2}||\tilde{\bxi}||^2$. The function $W$
is radially unbounded, positive definite, and satisfies $\dot{W}%&%=(1-\bth)\dot{V}_x+\bth\dot{V}_{\xi}\\
=(1-\lambda)\tilde{\mathbf{x}}^\top \dot{\tilde{\mathbf{x}}}+\lambda \tilde{\bxi}^\top\dot{\tilde{\bxi}}$. Let $\bh(\bar{\bth},\tilde{\mathbf{x}}):=\nabla_{\bx} f(\tilde{\mathbf{x}}+\mathbf{x}^*,\bar{\bth})$ and thus $\Bar{\bxi} = \Tilde{\bxi} + \bh(\bar{\bth},\tilde{\mathbf{x}})$. We rewrite the $\tilde{\mathbf{x}}$-error dynamics \eqref{eq:error:x} compactly as 
\begin{align} \label{eq:xtilde}  
\dot{\tilde{\mathbf{x}}}= \bbf_x(\tilde{\mathbf{x}},\tilde{\bxi}+\bh(\bar{\bth},\tilde{\mathbf{x}}))-\dot{\mathbf{x}}^*.
\end{align}
It follows that
\begin{align}
&\tilde{\mathbf{x}}^\top\dot{\tilde{\mathbf{x}}}=\tilde{\mathbf{x}}^\top (\bbf_x(\tilde{\mathbf{x}},\tilde{\bxi}+\bh(\bar{\bth},\tilde{\mathbf{x}}))-\dot{\mathbf{x}}^*) \nonumber \\
%&=\tilde{x}^\top (f_x(\tilde{x},\bxi)-\dot{x}^*\pm \bbf_x(\tilde{x},h(\bth,\tilde{x})))\\
=&\tilde{\mathbf{x}}^\top (\bbf_x(\tilde{\mathbf{x}},\tilde{\bxi}+\bh(\bar{\bth},\tilde{\mathbf{x}}))-\bbf_x(\tilde{\mathbf{x}},\bh(\bar{\bth},\tilde{\mathbf{x}}))) \nonumber \\
&~~~+\tilde{\mathbf{x}}^\top\left( \bbf_x(\tilde{\mathbf{x}},\bh(\bar{\bth},\tilde{\mathbf{x}}))-\dot{\mathbf{x}}^*\right)\nonumber  \\
\leq&||\tilde{\mathbf{x}}||\cdot||\dot{\mathbf{x}}^*||+||\tilde{\mathbf{x}}||\cdot||\bbf_x(\tilde{\mathbf{x}},\tilde{\bxi}+\bh(\bar{\bth},\tilde{\mathbf{x}}))-\bbf_x(\tilde{\mathbf{x}},\bh(\bar{\bth},\tilde{\mathbf{x}}))|| \nonumber   \\
&~~~+\tilde{\mathbf{x}}^\top \bbf_x(\tilde{\mathbf{x}},\bh(\bar{\bth},\tilde{\mathbf{x}})) \nonumber  \\
\leq & ||\tilde{\mathbf{x}}||\cdot||\dot{\mathbf{x}}^*||+k_x\alpha_x||\tilde{\mathbf{x}}||\cdot||\tilde{\bxi}||+\tilde{\mathbf{x}}^\top \bbf_x(\tilde{\mathbf{x}},\bh(\bar{\bth},\tilde{\mathbf{x}})), \label{eq:xxbound}
\end{align}
where we use the fact that $\mathcal{P}_{\mathcal{X}}(\cdot)$ is Lipschitz continuous with unitary Lipschitz constant \cite[Proposition 2.4.1]{clarke1990optimization}, such that:
\begin{align}\label{lipschitzproj}
\ ||\bbf_x(\tilde{\mathbf{x}},\bar{\bxi})-\bbf_x(\tilde{\mathbf{x}},\bh(\bar{\bth},\tilde{\mathbf{x}}))||&\leq \  k_x\alpha_x||\bar{\bxi}-h(\bar{\bth},\mathbf{x})||\nonumber\\
&=k_x\alpha_x||\tilde{\bxi}||.
\end{align} 
Using the uniform strong convexity and Lipschitz properties of Assumption \ref{ass:track} and the same steps of the proof of \cite[Theorem 4]{gao2003exponential}, we obtain 
\begin{equation*}
\tilde{\bx}^\top \bbf_x(\tilde{\mathbf{x}},\bh(\bar{\bth},\tilde{\mathbf{x}}))\leq -\tilde{k}_x||\tilde{\bx}||^2,
\end{equation*}
where $\tilde{k}_x=\alpha_x(\gamma-\frac{\alpha_x\ell^2}{4})$. Thus, for all $\alpha_x\in(0, \hat{\alpha}_x)$ with $\hat{\alpha}_x=\frac{4\gamma}{\ell^2}$ such that $\tilde{k}_x$ is positive, we obtain
\begin{align}\label{boundx1}
\tilde{\mathbf{x}}^\top \dot{\tilde{\mathbf{x}}}\leq -\tilde{k}_x||\tilde{\mathbf{x}}||^2+k_x\alpha_x||\tilde{\mathbf{x}}||\cdot||\tilde{\bxi}||+m||\tilde{\bx}||\cdot||\dot{\bar{\bth}}||,
\end{align}
where $m:=\max_{\bar{\bth}\in\Theta}||\nabla \bd(\bar{\bth})||$ and $\bd(\cdot)$ is defined in \eqref{eq:dtheta}. Note that $m<\infty$ because $\Theta$ is compact and $\bd$ is continuously differentiable by Assumption \ref{ass:track}. From \eqref{eq:error:xi}, we have 
\begin{equation*}
\tilde{\bxi}^\top \dot{\tilde{\bxi}}=-\frac{1}{\varepsilon_{\xi}}||\tilde{\bxi}||^2-\tilde{\bxi}^\top \frac{d}{dt}\nabla_{\bx} f(\tilde{\mathbf{x}}+\mathbf{x}^*,\bar{\bth}),
\end{equation*}
and note that
\begin{align*}
\frac{d}{dt}\nabla_{\bx}f(\tilde{\mathbf{x}}+\mathbf{x}^*,\bth)&=H_f(\bar{\bth},\tilde{\mathbf{x}}+\mathbf{x}^*)(\dot{\tilde{\mathbf{x}}}+\dot{\mathbf{x}}^*)\\&+\frac{\partial}{\partial\bth} \nabla_{\bx} f(\tilde{\mathbf{x}}+\mathbf{x}^*,\bar{\bth})\dot{\bar{\bth}},
\end{align*}
where $H_f(\cdot,\cdot)$ is the Hessian matrix of $f$. By \eqref{eq:globallipschitz1} in Assumption \ref{ass:track}, we have that $||H_f(\bar{\bth},\bar{\mathbf{x}})||\leq \ell$ for all $\bar{\mathbf{x}}\in\mathbb{R}^n$ and all $\bar{\bth}\in\Theta$, and since $\dot{\tilde{\mathbf{x}}}+\dot{\mathbf{x}}^*=\bbf_x(\tilde{\mathbf{x}},\tilde{\bxi}+\bh(\bar{\bth},\tilde{\mathbf{x}}))$ by \eqref{eq:xtilde}, we obtain
\begin{align}
&||H_f(\bar{\bth},\mathbf{\tilde{x}}+\mathbf{x}^*)(\dot{\tilde{\mathbf{x}}}+\dot{\mathbf{x}}^*)|| \nonumber\\
\leq &\ell ||\bbf_x(\tilde{\mathbf{x}},\bar{\bxi})-\bbf_x(\tilde{\mathbf{x}},\bh(\bar{\bth},\tilde{\mathbf{x}}))+\bbf_x(\tilde{\mathbf{x}},\bh(\bar{\bth},\tilde{\mathbf{x}}))||\notag\\
\leq  & \ell k_x\alpha_x||\tilde{\bxi}||\!+\!\ell||\bbf_x(\tilde{\mathbf{x}},\bh(\bar{\bth},\tilde{\mathbf{x}}))||
\leq \ell k_x\alpha_x||\tilde{\bxi}||\!+\!c\ell||\tilde{\bx}||,\label{usefullemma}
\end{align}
where $c:=k_x(2+\alpha_x\ell)$. The second inequality above is due to  \eqref{lipschitzproj}. The last inequality is because $\bbf_x(\bm{0},\bh(\bar{\bth},\bm{0}))=\bm{0}$  and
\begin{align*}
&\frac{1}{k_x} \cdot|| \bbf_x(\tilde{\mathbf{x}},\bh(\bar{\bth},\tilde{\mathbf{x}}))-\bbf_x(\bm{0},\bh(\bar{\bth},\bm{0}))|| \\
%&\Bigg|\Big(-x+\mathcal{P}_{\mathcal{X}}(x-\alpha \nabla f(x))\Big)-\Big(-x^*+\mathcal{P}_{\mathcal{X}}(x^*-\alpha \nabla f(x^*)) \Big)\Bigg|\\
\leq &||\tilde{\mathbf{x}}||\!+\!||\mathcal{P}_{\mathcal{X}}(\bar{\mathbf{x}}\!-\!\alpha_x \nabla_{\bx} f(\bar{\mathbf{x}},\bar{\bth}))\!-\!\mathcal{P}_{\mathcal{X}}(\mathbf{x}^*\!-\!\alpha_x \nabla_{\bx} f(\mathbf{x}^*,\bar{\bth}))| \\
\leq & ||\tilde{\mathbf{x}}||+||\tilde{\mathbf{x}}||+\alpha_x||\nabla_{\bx} f(\bar{\mathbf{x}},\bar{\bth})-\nabla_{\bx} f(\mathbf{x}^*,\bar{\bth})||\\
\leq  & 2||\tilde{\mathbf{x}}||+\alpha_x \ell||\tilde{\mathbf{x}}|| =  \frac{c}{k_x} ||\tilde{\bx}||. 
\end{align*}

Combining the above inequalities and using \eqref{boundpartial}, we obtain
\begin{equation}\label{xibound1}
\tilde{\bxi}^\top \dot{\tilde{\bxi}}\leq -\left(\frac{1}{\varepsilon_{\xi}}-\ell k_x\alpha_x\right)||\tilde{\bxi}||^2+c\ell||\tilde{\mathbf{x}}||\,||\tilde{\bxi}||+M||\tilde{\bxi}||\,||\dot{\bar{\bth}}||.
\end{equation}
Using \eqref{boundx1} and \eqref{xibound1}:
\begin{align}
\dot{W}&\leq (1\!-\!\lambda)k_x\alpha_x||\tilde{\mathbf{x}}||\,||\tilde{\bxi}|+m(1\!-\!\lambda)||\tilde{\bx}||\,||\dot{\bar{\bth}}||\!+\!M\lambda||\tilde{\bxi}||\,||\dot{\bar{\bth}}||\notag\\
&-\tilde{k}_x(1\!-\!\lambda)||\tilde{\mathbf{x}}||^2-\lambda\big(\frac{1}{\varepsilon_{\xi}}\!-\!\ell k_x\alpha_x\big)||\tilde{\bxi}||^2+c\ell\lambda||\tilde{\mathbf{x}}||\,||\tilde{\bxi}||\notag\\
&=-\tilde{\mathbf{y}}^\top Q \tilde{\mathbf{y}}+\bar{c}||\tilde{\mathbf{y}}||\,||\dot{\bar{\bth}}||,\label{Wbound}
\end{align}
where $\tilde{\mathbf{y}}=(||\tilde{\mathbf{x}}||,||\tilde{\bxi}||)$ and $Q$ is given by the $2\times2$ matrix
\begin{equation*}
Q=\left[\begin{array}{cc}
\tilde{k}_x(1-\lambda) & -\frac{1}{2}(1-\lambda)k_x\alpha_x-\frac{c\ell\lambda}{2}\\
-\frac{1}{2}(1-\lambda)k_x\alpha_x-\frac{c\ell\lambda}{2}  & \big(\frac{1}{\varepsilon_{\xi}}-\ell k_x\alpha_x\big)\lambda
\end{array}\right],
\end{equation*}
and $\bar{c}=\max\{m(1-\lambda),M\lambda\}$. This matrix is positive definite whenever
\begin{equation*}
\lambda(1-\lambda)\tilde{k}_x \left(\frac{1}{\varepsilon_{\xi}}-\ell k_x\alpha_x\right)>\frac{1}{4}\Big[(1-\lambda)k_x\alpha_x+\lambda c\ell\Big]^2,
\end{equation*}
which is equivalent to
\begin{equation*}
\frac{1}{\varepsilon_{\xi}}>\frac{0.25 \Big[(1-\lambda)k_x\alpha_x+\lambda c\ell\Big]^2}{\lambda(1-\lambda)\tilde{k}_x}+\ell k_x\alpha_x=:\frac{1}{\hat{\varepsilon}_{\xi}}.
\end{equation*}
This guarantees the existence of $\hat{\varepsilon}_{\xi}>0$ such that for all $\varepsilon_{\xi}\in(0,\hat{\varepsilon}_{\xi})$, the matrix $Q$ is positive definite, and, by \eqref{Wbound}, $W$ in \eqref{eq:W:lya} is a smooth input-to-state stability (ISS) Lyapunov function for the error dynamics \eqref{errordynamicsaux} with respect to the input $\dot{\bar{\bth}}$. This establishes the bound of Lemma \ref{auxlemma9}. Since $||\dot{\bth}||=\varepsilon_{\theta} ||\Pi(\bth)||\leq \varepsilon_{\theta} \sigma$, where $\sigma$ is such that $\Pi(\Theta)\subset\sigma\mathbb{B}$, the result also implies uniform ultimate boundedness (UUB) with ultimate bound proportional to $\varepsilon_{\theta}$.  \hfill $\blacksquare$

Since Lemma \ref{auxlemma9} directly establishes ISS for the nominal error averaged dynamics of \eqref{auxsysthm22} with $\mathcal{O}(\varepsilon_a)=0$, evolving in $\mathbf{C}_3$, the perturbed error average system \eqref{auxsysthm22} renders the origin semi-globally practically ISS as $\varepsilon_a\to0^+$. We can now directly link the stability properties of the average dynamics \eqref{auxsysthm22} and the original dynamics \eqref{eq:pzo} via standard averaging results for ISS systems \cite[Theorem 4]{wang2012input}. The practical safety property follows directly by Lemma \ref{lem:forward} since the dynamics \eqref{eq:pzo:x} are independent of $\bth$ (see the proof of Lemma \ref{lem:forward}).  

\vspace{0.1cm}\noindent
\subsection{Proof of Theorem \ref{thm:switchedprojected}}\label{sec:proof:switch}
For each constant $q$, the average dynamics of the P-GZO dynamics \eqref{eq:pzo} are given by
\begin{subequations}\label{averagedwitched}
\begin{align}
\dot{\bar{\mathbf{x}}}&=k_x \mathcal{P}_{\mathcal{X}}\left(\mathbf{\bar{x}}-\alpha_x \bar{\bxi}\right)-k_x\mathbf{\bar{x}},\\
\dot{\mathbf{\bar{\bxi}}}&=\frac{1}{\varepsilon_{\xi}}\left(-\mathbf{\bar{\bxi}}+\nabla f_q(\mathbf{\bar{x}})+\mathcal{O}(\varepsilon_a)\right),
\end{align}
\end{subequations}
which can be seen as a $\mathcal{O}(\varepsilon_a)$-perturbed two-time scale system with respect to the small parameter $\varepsilon_{\xi}$. The following lemma is needed to prove Theorem \ref{thm:switchedprojected}. We use the notation $\bs:=(\bx,\bxi)$ and $\mathcal{W}_1^*:=\{(\mathbf{x},\bxi)\in\mathbb{R}^{2n}:\mathbf{x}=\mathbf{x}^*,\,\bxi=\bxi^*:=\nabla f(\mathbf{x}^*)\}$.
\begin{lemma}\label{switchedsystem}
Consider the system  \eqref{averagedwitched} with $\mathcal{O}(\varepsilon_a)=0$. Then, there exist a function $V_q$ and $k_x,\hat{\varepsilon}_{\xi},c_1,c_2>0$ such that for all $\varepsilon_{\xi}\in(0,\hat{\varepsilon}_{\xi})$, we have
\begin{subequations}
\begin{align}
&c_1||\bs||^2_{\mathcal{W}_1^*}\leq V_q(\bs)\leq c_2||\bs||^2_{\mathcal{W}_1^*}\\
& \langle\nabla V_q(\bs),~\dot{\bs} \rangle\leq -\lambda_q V_q(\bs),
\end{align}
\end{subequations}
for all $\bs\in\mathcal{X}\times\mathbb{R}^n$ and all $q\in Q$. \QEDB
\end{lemma}
\textbf{Proof:} For each mode $q\in Q$, we consider the Lyapunov function $V_q(\bs)=\frac{1}{2}||\bx-\bx^*||^2+\frac{1}{2}||\bxi-\nabla f_q(\bx)||^2$. Since
$
||\bs||^2_{\mathcal{W}_1^*}=||\bx-\bx^*||^2+||\bxi-\nabla f_q(\bx^*)||^2$  
and  $\nabla f_q(\bx^*)=\bxi^*$,
\begin{align*}
&||\bxi-\nabla f_q(\bx)|| = ||\bxi-\nabla f_q(\bx^*)+\nabla f_q(\bx^*)-\nabla f_q(\bx) ||\\
\leq & ||\bxi-\bxi^*||+\ell_q ||\bx-\bx^*||,
\end{align*}
where we used the Lipschitz property of $\nabla f_q$. Hence, there exists $\bar{c}_q>0$ such that $V_q(\bs)\leq \bar{c}_q||\bs||^2_{\mathcal{W}_1^*}$ for all $\bs\in\mathcal{X}\times\mathbb{R}^n$. Similarly,
\begin{align*}
||\bxi-\nabla f_q(\bx^*) ||&\leq ||\bxi-\nabla f_q(\bx)||+||\nabla f_q(\bx)-\nabla f_q(\bx^*)||\\
&\leq ||\bxi-\nabla f_q(\bx)||+\ell_q||\bx-\bx^*||.
\end{align*}
Therefore,
\begin{align*}
||\bxi-\nabla f_q(\bx^*)||^2&\leq \left(||\bxi-\nabla f_q(\bx)||+\ell_q||\bx-\bx^*||\right)^2\\
&\leq 2||\bxi-\nabla f_q(\bx)||^2+2\ell_q^2||
\bx-\bx^*||^2.
\end{align*}
Adding $||\bx-\bx^*||^2$ to both sides and dividing by $d_q=2(\ell_q^2+1)$ leads to
\begin{equation}
\frac{1}{2d_q}||\bs||^2_{\mathcal{W}_1^*}\leq V_q(\bs),~~\forall~ \bs\in\mathcal{X}\times\mathbb{R}^n.
\end{equation}
It follows that $c_2=\max_{q\in Q}\bar{c}_q$, and $c_1=1/(2\max_{q\in Q}{d_q})$. Next, note that
\begin{equation*}
\dot{V}_q\!=\!(\bx-\bx^*)^\top\dot{\bx}-(\bxi-\nabla f_q(\bx))^\top H_{f_q}(\bx)\dot{\bx}+(\bxi-\nabla f_q(\bx))^\top \dot{\bxi}. 
\end{equation*}
The first term of $\dot{V}_q$ satisfies
\begin{align*}
    & (\bx-\bx^*)^\top \dot{\bx}=(\bx-\bx^*)^\top (\dot{\bx}- \dot{\bx}_r + \dot{\bx}_r)\\
    \leq\, &(\bx-\bx^*)^\top \dot{\bx}_r+||\bx-\bx^*||\,||\dot{\bx}-\dot{\bx}_r||\\
    \leq\, &-\tilde{k}_x||\bx-\bx^*||^2+k_x\alpha_x||\bx-\bx^*||\,||\bxi-\nabla f_q(\bx)||,
\end{align*}
where $\dot{\bx}_r$ denotes the right-hand side of \eqref{eq:pzo:x} with $\bxi=\nabla f_q(\bx)$ and   we used \cite[Theorem 4]{gao2003exponential} and \eqref{lipschitzproj} to obtain the last inequality. Next, using \eqref{usefullemma}, the second term of $\dot{V}_q$ can be bounded as
\begin{align*}
-(\bxi-\nabla f_q(\bx))^\top H_{f_q}(\bx)\dot{\bx}&\leq \ell k_x\alpha_x||\bxi-\nabla f_q(\bx)||^2\\
&~~+c\ell||\bxi-\nabla f_q(\bx)||\,||\bx-\bx^*||.
\end{align*}
The last term of $\dot{V}_q$ satisfies 
$(\bxi-\nabla f_q(\bx))^\top \dot{\bxi}=-\frac{1}{\varepsilon_{\xi}}||\bxi-\nabla f_q(\bx)||^2$. Therefore, $\dot{V}_q$ satisfies the same upper bound of \eqref{Wbound} with $\dot{\bth}=0$,  and there exists $\hat{\varepsilon}_{\xi}>0$ sufficiently small such that for all ${\varepsilon}_{\xi}\in(0,\hat{\varepsilon}_{\xi})$, $Q$ is positive definite and 
\begin{equation}
\dot{V}_q\leq -k V_q(\bs),~\quad \forall\bs\in\mathcal{X}\times\mathbb{R}^n.
\end{equation}
This establishes the result of  Lemma \ref{switchedsystem}. \hfill $\blacksquare$

The result of Lemma \ref{switchedsystem}, in conjunction with \cite[Exercise 3.22]{Goebel:12}, guarantees the existence of a $\tau_d$ sufficiently large such that, the hybrid dynamical system with flows 
\begin{subequations}\label{averagehybrid}
\begin{align}
\dot{\bar{\mathbf{x}}}&=k_x \mathcal{P}_{\mathcal{X}}\left(\mathbf{\bar{x}}-\alpha_x \bar{\bxi}\right)-k_x\mathbf{\bar{x}},\\
\dot{\mathbf{\bar{\bxi}}}&=\frac{1}{\varepsilon_{\xi}}\left(-\mathbf{\bar{\bxi}}+\nabla f_{\bar{q}}(\mathbf{\bar{x}})\right)\\
\dot{\bar{q}}&=0\\
\dot{\bar{\tau}}&\in\left[0,\frac{1}{\tau_d}\right]
\end{align}
\end{subequations}
evolving on the flow set $C=(\mathcal{X}\times\mathbb{R}^n)\times [0,N_0]\times\mathcal{Q}$, and jumps
\begin{equation}
\bar{\mathbf{x}}^+=\mathbf{x},~~\bar{\mathbf{\xi}}^+=\bar{\mathbf{\xi}},~~\bar{q}^+\in Q,~~\bar{\tau}^+=\bar{\tau}-1.
\end{equation}
evolving on the jump set $D=(\mathcal{X}\times\mathbb{R}^n)\times [1,N_0]\times\mathcal{Q}$, renders the set $\mathcal{W}_1^*\times[0,N_0]\times\mathcal{Q}$ UGAS. Note that this hybrid system is well-posed in the sense of \cite[Def. 6.29]{Goebel:12}, as it satisfies the hybrid basic conditions \cite[Assumption 6.5]{Goebel:12}. Moreover, by \cite[Prop. 6.10]{Goebel:12}, the existence of solutions (as hybrid arcs \cite[Def. 2.4]{Goebel:12}) from all initial conditions in $C\cup D$ is guaranteed.

In turn, by robustness properties of well-posed hybrid systems, the $\mathcal{O}(\varepsilon_a)$-perturbation of this nominal average system renders the same set SGPAS as $\varepsilon_a\to0^+$ \cite[Theorem 7.21]{Goebel:12}. Therefore, the result of Theorem \ref{thm:switchedprojected} follows now directly by an application of averaging theory for perturbed hybrid systems \cite[Theorem 7]{PovedaNaLi2019}.  
\subsection{Proof of Theorem \ref{thm:dpzo}} \label{sec:proof:dpzo}
Since the DP-GZO dynamics \eqref{eq:dpzo} is a discontinuous ODE, we consider its  Krasovskii regularization defined in \eqref{krasovskiiregularization}, which only affects the right-hand side of $\dot{\mathbf{x}}$:
\begin{equation}\label{dproof1}
\dot{\mathbf{x}}\in K(\mathbf{z}),\quad\dot{\bxi}=\frac{1}{\varepsilon_\xi} \big(-\bxi+\frac{2}{\epsilon_a}f(\hat{\mathbf{x}})\hat{\bmu}\big),
\end{equation}
where $\mathbf{z}:=(\mathbf{x},\bxi)$. Since $\mathcal{X}$ is closed and convex, and $f$ is continuously differentiable, by \cite[Theorem 4.2]{thesisadrian}, every solution of \eqref{dproof1} is also a solution of the DP-GZO dynamics \eqref{eq:dpzo}, and vice versa. Moreover, since the dynamics of $\mathbf{x}$ are  independent of $\bmu$, system \eqref{dproof1} is in standard form for the application of averaging theory \cite[Definition 7]{PovedaTAC17B}. In particular, similar to Lemma \ref{lemma:aveg2}, we compute the average dynamics of \eqref{dproof1} along $t\to\bmu(t)$ and obtain
\begin{align} \label{dproof2}
\dot{\mathbf{\bar{\bx}}}\in K(\mathbf{\bar{\bz}}),\quad \dot{\bar{\bxi}}=\frac{1}{\varepsilon_{\xi}} \left(-\bar{\bxi}+\nabla f(\mathbf{\bar{x}})+\mathcal{O}(\varepsilon_a)\right).
\end{align}
which can be seen as an $\mathcal{O}(\varepsilon_a)$-perturbed two-time scale set-valued dynamical system. We will first study the stability properties of this system by analyzing the nominal unperturbed dynamics corresponding to $\mathcal{O}(\varepsilon_a)=0$.
\begin{lemma}\label{auxlemmaD}
Under the assumptions of Theorem \ref{thm:dpzo}, system \eqref{dproof2} with $\mathcal{O}(\varepsilon_a)=0$ and flow set $ \mathcal{X}\times\mathbb{R}^n$ renders the point $\bz^*=(\mathbf{x}^*,\nabla f(\mathbf{x}^*))$ UGAS. \QEDB
\end{lemma}

\vspace{0.1cm}\noindent 
\textbf{Proof}: Using the equivalence between Krasovskii and Caratheodory solutions for well-posed projected gradient systems \cite[Theorem 4.2]{thesisadrian}, we consider the dynamics
\begin{equation}\label{dproof4}
\dot{\mathbf{\bar{\bx}}}=k_x\mathcal{P}_{T_{\mathcal{X}}(\bar{\bx})}(-\bar{\bxi}),~~\dot{\bar{\bxi}}=\frac{1}{\varepsilon_\xi} \left(-\bar{\bxi}+\nabla f(\mathbf{\bar{\bx}})\right),
\end{equation}
and the composite Lyapunov function with $\lambda\in(0,1)$
\begin{align}
    V(\bar{\mathbf{z}})=(1-\lambda)(f(\bar{\mathbf{x}})-f(\bx^*))+\lambda\frac{1}{2}||\bar{\bxi}-\nabla f(\bar{\mathbf{x}})||^2,
\end{align}
 which is continuously differentiable, radially unbounded, and positive definite with respect to $\mathbf{z}^*$ in $\sX\times\R^n$.
 
 We proceed to compute the inner product $\langle \nabla V, \dot{\bar{\mathbf{z}}} \rangle$, where $\dot{\bar{\mathbf{z}}}=(\dot{\bar{\mathbf{x}}},\dot{\bar{\bxi}})$. To do this, we use the fact that for any regular set $\mathcal{X}$, and any ${\mathbf{x}}\in\mathcal{X}$, $\bm{\nu}\in\mathbb{R}^n$, there exists a unique $\bm{\eta}\in N_{\mathcal{X}}(\bx)$ such that $\mathcal{P}_{T_{\mathcal{X}}(\bx)}(\bm{\nu})=\bm{\nu}-\bm{\eta}$, $\bm{\eta}^\top (\bm{\nu}-\bm{\eta})=0$, and $\bm{\nu}^\top (\bm{\nu}-\bm{\eta})=||\bm{\nu}-\bm{\eta}||^2$, 
 \cite[Lemma C.3]{hauswirth2020timescale}.
 Thus using $\bm{\nu}=-\bar{\bxi}$, $\tilde{\mathbf{v}}(\bar{\mathbf{x}}):=\bm{\nu}-\bm{\eta}$, and $\bh(\bar{\mathbf{x}})=\bar{\bxi}-\nabla f(\bar{\mathbf{x}})$ we obtain:
\begin{align*}
& \langle \nabla V, \dot{\bar{\mathbf{z}}} \rangle\\
 = & k_x(1\!-\!\lambda)\nabla f(\bar{\mathbf{x}})^\top \tilde{\bv}(\bar{\mathbf{x}})\!-\!k_x\lambda \bh(\bar{\mathbf{x}})^\top\! H_f(\bar{\mathbf{x}}) \tilde{\bv}(\bar{\mathbf{x}}) \!+\!\lambda \bh(\bar{\mathbf{x}})^\top \dot{\bxi}\\
=& k_x(1\!-\!\lambda)\nabla f(\bar{\mathbf{x}})^\top \tilde{\bv}(\bar{\mathbf{x}})\!-\!k_x\lambda \bh(\bar{\mathbf{x}})^\top\! H_f(\bar{\mathbf{x}}) \tilde{\bv}(\bar{\mathbf{x}}) \!-\!\frac{\lambda}{\varepsilon_{\xi}} ||\bh(\bar{\mathbf{x}})||^2.
\end{align*}
To upper-bound the first term, we note that
\begin{align*}
&\nabla f(\bar{\mathbf{x}})^\top \tilde{\bv}(\bar{\mathbf{x}})
=(\nabla f(\bar{\mathbf{x}})-\bar{\bxi})^\top \tilde{\bv}(\bar{\mathbf{x}})+\bar{\bxi}^\top \tilde{\bv}(\bar{\mathbf{x}})\\
%&\leq ||\bh(\bar{\mathbf{x}})||\,||\tilde{\bv}(\bar{\mathbf{x}})||+\bar{\bxi}^\top (\bm{\nu}-\bm{\eta})\\
\leq& ||\bh(\bar{\mathbf{x}})||\,||\tilde{\bv}(\bar{\mathbf{x}})||-\bm{\nu}^\top (\bm{\nu}-\bm{\eta})\\
\leq & ||\bh(\bar{\mathbf{x}})||\,||\tilde{\bv}(\bar{\mathbf{x}})||-||\bm{\nu}\!-\!\bm{\eta}||^2 = ||\bh(\bar{\mathbf{x}})||\,||\tilde{\bv}(\bar{\mathbf{x}})|| -||\tilde{\bv}(\bar{\mathbf{x}})||^2.
\end{align*}
%
%
%
%
% \begin{align*}
% \nabla f(\mathbf{x})^\top \tilde{v}(\mathbf{x}) \leq k_x L|\xi-\nabla f(\mathbf{x})||\tilde{v}(\mathbf{x})|-k_x |\tilde{v}(\mathbf{x})|^2
% \end{align*}
% %
% where 
% %
% \begin{equation}
% \tilde{v}(x)= v-\eta,~v=-\xi,~\eta\in N_x\mathcal{X}.
% \end{equation}
%
Moreover, since by assumption $\nabla f$ is $\ell$-globally Lipschitz, the second term of $\dot{V}$ satisfies 
$$(\bar{\bxi}-\nabla f(\bar{\mathbf{x}}))^\top H_f(\bar{\mathbf{x}}) \tilde{\bv}(\bar{\mathbf{x}}) \leq \ell||\bh(\bar{\mathbf{x}})||\,||\tilde{\bv}(\bar{\mathbf{x}})||.$$ Therefore, defining $\tilde{\bq}(\bar{\mathbf{z}}):=(\tilde{\bv}(\bar{\mathbf{x}}),\bh(\bar{\mathbf{x}}))$, we obtain:
\begin{equation}\label{Lyapunovdiscon}
 \langle \nabla V(\bar{\mathbf{z}}), \dot{\bar{\mathbf{z}}} \rangle \leq - \tilde{\bq}(\bar{\mathbf{z}})^\top Q\, \tilde{\mathbf{q}}(\bar{\mathbf{z}}),
\end{equation}
where
\begin{equation*}
Q=\left[\begin{array}{cc}
k_x(1-\lambda) & -\frac{1}{2}(k_x\lambda \ell+k_x(1-\lambda)) \\
-\frac{1}{2}(k_x\lambda 
\ell+k_x(1-\lambda)) & \lambda \frac{1}{\varepsilon_{\xi}}
\end{array}\right].
\end{equation*}
This matrix is positive definite whenever $\lambda (1-\lambda) \frac{k_x}{\varepsilon_{\xi}}>\frac{1}{4}[k_x\lambda\ell +k_x(1-\lambda)]^2$, which can be satisfied for sufficiently small values of $\varepsilon_{\xi}$. Since $\tilde{\bq}(\mathbf{z})=0$ if and only if $\mathbf{z}=\mathbf{z}^*=(\mathbf{x}^*,\bxi^*)$, we obtain that $\sX^*$ is uniformly globally asymptotically stable (UGAS) for \eqref{dproof4}. %By equivalence between Krasovskii and Caratheodory solutions, the stability and convergence properties are uniform.
\hfill $\blacksquare$

\vspace{0.1cm}
By equivalence between Krasovskii and Caratheodory solutions, the result of Lemma \ref{auxlemmaD} guarantees UGAS for the Krasovskii regularization of \eqref{dproof4}, which is precisely \eqref{dproof2} with $\mathcal{O}(\varepsilon_a)=0$. Since, by construction, this system is well-posed (outer-semi-continuous, locally bounded, and convex-valued), the set $\sX^*$ is semi-globally practically asymptotically stable for \eqref{dproof2} as $\varepsilon_a\to0^+$. The stability result of Theorem \ref{thm:dpzo} follows now by a direct application of averaging for non-smooth systems of the form \eqref{dproof1} \cite[Lemma 6]{PovedaTAC17B}. \hfill $\blacksquare$

\section{CONCLUSION}\label{sec:conclusion}
In this paper, we introduce a class of continuous-time projected zeroth-order (P-ZO) dynamic methods for solving generic constrained optimization problems with both hard and asymptotic constraints. In these problems, the mathematical forms of the objective and constraint functions are unknown, and only their function evaluations are available. Consequently, the proposed P-ZO methods can be interpreted as model-free feedback controllers that guide a black-box plant toward optimal steady states defined by an optimization problem using only measurement feedback. We consider both continuous and discontinuous projection maps, establishing the stability and robustness of the proposed P-ZO methods. Additionally, we analyze their dynamic tracking performance under time-varying settings and switching cost functions. Future research directions include the study of non-convex and switching cost functions with no-common critical points, problems with closed rather than compact sets of saddle points, projected exploration dithers, and the practical implementation of the P-ZO methods in practical problems where the cost to be minimized is the output of a dynamic plant.

%we propose the P-PDZD method to solve the generic constrained optimization problems  with hard and asymptotic constraints in a model-free feedback manner.  Using only zeroth-order feedback, the proposed method can be interpreted as the model-free feedback controller that autonomously drives a black-box system to the solution of the optimization problem. We prove the semi-global practical asymptotic stability and structural robustness of the P-PDZD and present the decentralized version of P-PDZD when applied to multi-agent problems. The numerical simulations on the optimal voltage control problem with square-wave probing signals demonstrate the optimality, robustness, and dynamic tracking capability of the P-PDZD. For future work, we will incorporate the plant dynamics and study the discrete-time implementation of the P-PDZD.

%\newpage

\appendices

\section{AUXILIARY LEMMAS}

\label{sectionproofs}

\subsection{Proof of Lemma \ref{lem:forward}}
\label{app:lem:forward:pf}

%\textbf{Proof:} 
For the purpose of analysis, we take the flow set of \eqref{eq:pzo} to be $\mathbb{R}^n\times\mathbb{R}^n\times\mathbb{T}^n$, since otherwise there is nothing to prove. First, we let
$\bmu(0)\in \mathbb{T}^n$ and   $\mu_i(0)^2\!+\!\mu_{i+1}(0)^2\!=\!1$ for all $i\in\{1,3,\ldots,2n-1\}$. Since
\begin{align*}
\frac{d}{dt}\left(\mu_i(t)^2+\mu_{i+1}(t)^2\right)&=2\mu_i\dot{\mu}_i+2\mu_{i+1}\dot{\mu}_{i+1}\\
&=2 (\mu_i,\mu_{i+1})^\top \Lambda_i(\mu_i,\mu_{i+1})=0, 
\end{align*}
$\mathbb{T}^n$ is forward invariant for $\bmu(t)$ under \eqref{eq:pzo:mu}. Then, following the ideas of \cite[Theorem 3.2]{xia2000stability}, we define $\Phi(\mathbf{x}):=|\mathbf{x}-\mathcal{P}_{\mathcal{X}}(\mathbf{x})|^2$ and have
\begin{align*}
\dot{\Phi}&=2(\mathbf{x}-\mathcal{P}_{\mathcal{X}}(\mathbf{x}))^\top \dot{\mathbf{x}}\\
% &=k_x(\mathbf{x}-\mathcal{P}_{\mathcal{X}}(\mathbf{x}))^\top\left(\mathcal{P}_{\mathcal{X}}(\mathbf{x}-\alpha_x\mathbf{\xi})-\mathbf{x}\right)\\
&=2k_x(\mathbf{x}-\mathcal{P}_{\mathcal{X}}(\mathbf{x}))^\top\left(\mathcal{P}_{\mathcal{X}}(\mathbf{x}-\alpha_x\bm{\xi})-\mathbf{x}\right)\\
&=-2k_x(\mathbf{x}-\mathcal{P}_{\mathcal{X}}(\mathbf{x}))^\top\left(\mathbf{x}- \mathcal{P}_{\mathcal{X}}(\mathbf{x})\right)\\
&~~~-2k_x(\mathbf{x}-\mathcal{P}_{\mathcal{X}}(\mathbf{x}))^\top\left(\mathcal{P}_{\mathcal{X}}(\mathbf{x})-\mathcal{P}_{\mathcal{X}}(\mathbf{x}-\alpha_x{\bxi})\right)\\
&\leq  -2k_x |\mathbf{x}-\mathcal{P}_{\mathcal{X}}(\mathbf{x})|^2 = -2k_x\Phi(\mathbf{x}), 
\end{align*}
for all $\mathbf{x}\in\mathbb{R}^n$, where the first equality follows by \cite[Prop. 3.1]{thesisadrian}, and the inequality in the last step used the property  
 $(\mathbf{u}-P_{\mathcal{X}}(\mathbf{u}))^\top(P_{\mathcal{X}}(\mathbf{u})-\mathbf{v})\geq 0$ for all $\mathbf{u}\in\mathbb{R}^n$ and all $\mathbf{v}\in\mathcal{X}$.
 This implies that $\dot{\Phi}(\mathbf{x}(t))\leq -2k_x\Phi(\mathbf{x}(t))\leq 0$, for all $t\in\text{dom}(\mathbf{z})$. To show that $\bx(0)\in \sX$ implies  $\bx(t)\in \sX$ for all $t\in \text{dom}(\bz)$, suppose by contradiction that there exists  $t_2>t_1$ with   $t_2,t_1\in\text{dom}(\mathbf{z})$ such that    $\mathbf{x}(t)\in\mathcal{X}$ for all $t\in[0,t_1]$ and  $\mathbf{x}(t)\notin\mathcal{X}$ for all $t\in(t_1,t_2]$. Then, it follows that $\Phi(\mathbf{x}(t_1))=0$ and $\Phi(\mathbf{x}(t_2))>0$. But the mean value theorem implies the existence of a   $ \bar{t}\in(t_1,t_2)$ such that $\dot{\Phi}(\bar{t})=\frac{\Phi(t_2)-\Phi(t_1)}{t_2-t_1}>0$, which is a contradiction. Therefore, we conclude that if $\mathbf{z}(0)\in \mathbf{C}_1$, then $\mathbf{z}(t)\in \mathbf{C}_1$ for all $t\in\text{dom}(\mathbf{z})$. Since the input $\hat{\mathbf{x}}$ is defined via \eqref{eq:hatx} and $|\mu_i(t)|\leq 1$ for all $i$ and $t\geq0$, then $\hat{\mathbf{x}}(t)\in \mathcal{X}+\varepsilon_a\mathbb{B}$ for all $t\in\text{dom}(\mathbf{z})$. % \hfill $\blacksquare$

\subsection{Proof of Lemma \ref{lemma:aveg2}} \label{sec:app:lem}

%Based on the concrete formulation \eqref{eq:g2} of $\bg_2$,
First, consider the  integration on the first part of $\bq_2(y,{\bmu}(t))$. By the Taylor expansion of $f(\cdot)$, we have  $(\forall i\in[n])$ 
\begin{align*}
 & \frac{1}{T}\int_{0}^T \frac{2}{\varepsilon_a }  f\left(\bx+ \varepsilon_a \hat{\bmu}(t)\right) \hat{\bmu}_i(t)\, dt \\
 = &  \frac{1}{T}\int_{0}^T\! 
\frac{2}{\varepsilon_a }  \big[f(\bx)\!+\! \varepsilon_a  \nabla f(\bx)^\top \hat{\bmu}(t)
\! +\!  \sO(\varepsilon_a^2)\big] \hat{\bmu}_i(t) dt \\
= & \frac{1}{T}\int_{0}^T {2} \sum_{j=1}^n \big[\frac{\partial f(\bx)}{\partial x_j} \hat{\bmu}_j(t)\hat{\bmu}_i(t)\big]\, dt +\sO(\varepsilon_a)\\
= & \frac{\partial f(\bx)}{\partial x_i} \frac{\eta_d}{ T}\int_{0}^T \!   \hat{\bmu}_i(t)^2 \, dt +\sO(\varepsilon_a)= \frac{\partial f(\bx)}{\partial x_i}  +\sO(\varepsilon_a).
\end{align*}
Similarly, we have $(\forall j\in[m], i\in[n])$
\begin{align*}
  \frac{1}{T}\!\!\int_{0}^T\! \!\frac{2}{\varepsilon_a }   \lambda_jg_j(\bx\!+\! \varepsilon_a \hat{\bmu}(t)) \hat{\bmu}_i(t) dt\!= \! \lambda_j \frac{\partial g_j(\bx)}{\partial x_i} \!+\!\sO(\varepsilon_a).
\end{align*}
%{\color{red} there seems a problem here, the perturbation term should be $\sO(\lambda_j \epsilon_a)$?? So can not put the gradient of f and g together?}

As for the integration on the second part of $\bq_2(y,\mathbf{\mu}(t))$, i.e., $\bg(\hat{\bx}(t))$, each component of this integration is ($\forall j\in[m]$)
\begin{align*}
     &\, \frac{1}{T}\int_{0}^T \!g_j(\bx+ \epsilon_a\hat{\bmu}(t)) \, dt \\
= & \,\frac{1}{T}\int_{0}^T \!
g_j(\bx) + \epsilon_a  \nabla g_j(\bx)^\top \hat{\bmu}(t) +\sO(\varepsilon_a^2)\, dt= \,g_j(\bx) + \sO(\varepsilon_a^2). 
\end{align*}
Combining these two parts, Lemma \ref{lemma:aveg2} is proved.

\subsection{Proof of Lemma \ref{lemma:omega}}\label{app:omega}

%{\color{red} not finish yet. $\nabla J$ is not necessarily converges to 0 due to the projection.}

Take $\varepsilon_a$ sufficiently small such that $\sO(\epsilon_a)< 1$ in \eqref{eq:realave}. 
For each $\nu\in(0,1)$, there exists a time $T_1>0$ such that for any $t\geq T_1$,  $\beta(\Delta, t)\leq \frac{\nu}{4}$. Such  $T_1$ always exists because $\beta$ is a class-$\mathcal{KL}$ function, and thus  $|\bar{\by}(t)|_{\mathcal{Y}^*}\leq \frac{\nu}{2}$ for $t\geq T_1$ by \eqref{eq:y1con}. In addition, by the exponential input-to-state stability of the linear fast dynamics \eqref{eq:realave}, there exists $T_2>0$ such that for any $t\geq T_2$, every solution of \eqref{eq:realave} with $\bar{\bs}(0)\in [(\mathcal{Y}^*\!+\!\Delta \mathbb{B})\cap \mathcal{Y}]\!\times\! \Delta \B$ satisfies  $|\bar{\bxi}(t)|\leq \frac{\nu}{2} + \sup_{\tau\geq t_0} ||\mathbf{\ell}(\bar{\by}(\tau)) +\sO(\varepsilon_a)||\leq \frac{\nu}{2}+ M_2$. Thus, for all $t\geq \max\{T_1,T_2\}$, the trajectory $\bar{\bs}$ converges to a $\frac{\nu}{2}$-neighborhood of $\mathcal{Y}^*\times M_2\B$. Therefore, the Omega-limit set from $\sF\times M_3\B$ is nonempty and satisfies $\Omega(\sF\times M_3\B)\subset (\mathcal{Y}^*\times M_2\B) + \frac{\nu}{2} \B \subset \mathrm{int}(\sF\times M_3\B)  $.  By \cite[Corollary 7.7]{Goebel:12}, the set $\Omega(\sF\times M_3\B)$ is uniformly globally asymptotically stable for the average system \eqref{eq:realave} restricted to $\sF\times M_3\B$. %Thus Lemma \ref{lemma:omega} is proved.

\bibliography{IEEEabrv, ref}

\newpage
% \vspace{-1.2cm}
% \begin{IEEEbiography}[{\includegraphics[width=1in,height=1.25in,clip,keepaspectratio]{Xin.jpg}}]{Xin Chen} 
\begin{IEEEbiographynophoto}{Xin Chen}
is an Assistant Professor in the Department of Electrical and Computer Engineering at Texas A\&M University (TAMU). Prior to joining TAMU, he was a Postdoctoral Associate affiliated with MIT Energy Initiative at Massachusetts Institute of Technology. He received the Ph.D. degree in electrical engineering from Harvard University, the master’s degree in electrical engineering and two bachelor’s degrees in engineering and economics from Tsinghua University. Dr. Chen is a recipient of the IEEE PES Outstanding Doctoral Dissertation, IEEE Transactions on Smart Grid Top-5 Outstanding Papers, the Best Research Award at the 2023 IEEE PES Grid Edge Conference, the Outstanding Student Paper Award at the 2021 IEEE Conference on Decision and Control,  the Best Student Paper Award Finalist at the 2018 IEEE Conference on Control Technology and Applications, and the Best Conference Paper Award at the 2016 IEEE PES General Meeting.
% \end{IEEEbiography}
\end{IEEEbiographynophoto}

% \vspace{-1.2cm}
% \begin{IEEEbiography} [{\includegraphics[width=1in,height=1.25in,clip,keepaspectratio]{Poveda.jpg}}]{Jorge I. Poveda}
\begin{IEEEbiographynophoto}{Jorge I. Poveda}
is an Assistant Professor in the Department of Electrical and Computer Engineering at the University of California, San Diego. He received his M.Sc. and Ph.D. degrees in Electrical and Computer Engineering from UC Santa Barbara in 2016 and 2018, respectively, and was Postdoctoral Fellow at Harvard University during part of 2018. 
He has received the Donald P. Eckman Award, AFOSR Young Investigator Award, the NSF CRII and CAREER awards, the CCDC Best Ph.D. Thesis award and Outstanding Scholar Fellowship from UC Santa Barbara, and the 2023 IEEE Transactions on Control of Network Systems Best Paper Award. He has served as Associate Editor
for Automatica, IEEE Control System Letters, and for Nonlinear Analysis: Hybrid Systems.
% \end{IEEEbiography}
\end{IEEEbiographynophoto}
%%

% \vspace{-1.2cm}
% \begin{IEEEbiography}
% [{\includegraphics[width=1in,height=1.25in,clip,keepaspectratio]{Lina.jpg}}]{Na Li} 

\begin{IEEEbiographynophoto}{Na Li}
is the Gordon McKay Professor of Electrical Engineering and Applied Mathematics in the School of Engineering and Applied Sciences at Harvard University. She received her PhD degree in Control and Dynamical systems from the California Institute of Technology in 2013. In 2014, she was a postdoctoral associate of the Laboratory for Information and Decision Systems
at Massachusetts Institute of Technology.  She has received the Donald P. Eckman Award, ONR and AFOSR  Young Investigator Awards, NSF CAREER Award, the Harvard PSE Accelerator Award, and she was also a Best Student Paper Award finalist in the 2011 IEEE Conference on Decision and Control. She has served as Associate Editor for IEEE Transactions on Automatic Control, Systems and Control Letters, and IEEE Control System Letters.
%%
% \end{IEEEbiography}
\end{IEEEbiographynophoto}

% input{Old}

\end{document}